\documentclass[12pt]{article}
\usepackage{amsfonts}
\usepackage{hyperref}
\usepackage{amsmath,amssymb,amsthm}
\usepackage[square,numbers]{natbib}
\usepackage{pst-node,pst-tree}
\makeatletter
\newcommand\figcaption{\def\@captype{figure}\caption}
\makeatother

\newtheorem{theorem}{Theorem}

\newtheorem{corollary}[theorem]{Corollary}
\newtheorem{lemma}[theorem]{Lemma}

\theoremstyle{definition}

\newtheorem{example}[theorem]{Example}

\newtheorem{proposition}[theorem]{Proposition}
\newtheorem{remark}[theorem]{Remark}

\begin{document}

\title{Measures of the full Hausdorff dimension for a general Sierpi\'{n}ski carpet}
\author{Jung-Chao Ban\thanks{%
The first author is partially supported by the National Science Council, ROC
(Contract No NSC 100-2115-M-259-009-MY2).}\thanks{%
Corresponding author} \\
Department of Applied Mathematics,\\
National Dong Hwa University,\\
Hualien, 97401, Taiwan, R.O.C. \and Chih-Hung Chang\thanks{%
The second author is grateful for the partial support of the
National Science Council, ROC (Contract No NSC 100-2115-M-035-003).} \\
Department of Applied Mathematics, \\
Feng Chia University,\\
Taichung 40724, Taiwan, R.O.C. \and Ting-Ju Chen \\
Department of Applied Mathematics, \\
National Dong Hwa University, \\
Hualien, 97401, Taiwan, R.O.C.}
\maketitle

\begin{abstract}
The measure of the full dimension for a general Sierpi\'{n}ski carpet is
studied. In the first part of this study, we give a criterion for
the measure of the full Hausdorff dimension of a Sierpi\'{n}ski carpet.
Meanwhile, it is the conditional equilibrium measure of zero
potential with respect to some Gibbs measure $\nu _{\alpha }$ of
matrix-valued potential $\alpha\mathbf{N}$ (defined later). On one hand, this
investigation extends the result of \cite{Olivier2009} without
condition \textbf{(H)}. On the other hand, it provides a checkable
condition to ensure the existence and uniqueness of  the measure of the full
Hausdorff dimension for a general Sierpi\'{n}ski carpet.

In the second part of this paper we give a criterion for the  Markov projection
measure and estimate its number of steps by means of the induced
matrix-valued potential. The results enable us to answer some
questions which arise from \cite{Chazottes2003} and \cite{BarralFeng2009}
on the projection measure and factors.
\end{abstract}

\begin{flushleft}

\end{flushleft}

\section{Introduction and main results}

Let $\mathbb{T}^{2}=\mathbb{R}^{2}/\mathbb{Z}^{2}$ which is invariant under
the endomorphism%
\begin{equation*}
\mathbb{T}=\left(
\begin{array}{cc}
\mathbf{m} & 0 \\
0 & \mathbf{n}%
\end{array}%
\right)
\end{equation*}%
and corresponds to a shift of finite type. Denote the set of digits as follows:
\begin{equation}
D=\left\{ 1,\ldots ,\mathbf{m}\right\} \times \left\{ 1,\ldots ,\mathbf{n}%
\right\}   \label{32}
\end{equation}%
For $\left( d_{k}\right) _{k=1}^{\infty }\in D^{\mathbb{N}}$ there
corresponds a point in $\mathbb{T}^{2}$ via what may be called "\emph{%
base }$\mathbb{T}$\emph{\ representation}" \cite{Kenyon1996}.%
\begin{equation*}
R_{\mathbb{T}}\left( \left( d_{k}\right) \right) =\sum_{k=1}^{\infty }\left(
\begin{array}{cc}
\mathbf{m}^{-k} & 0 \\
0 & \mathbf{n}^{-k}%
\end{array}%
\right) d_{k}.
\end{equation*}%
Any $0$-$1$ matrix $A$ with rows and columns indexed by $D$ defines a shift
of finite type (SFT, for short) and let $K_{\mathbb{T}}(A)$ be its \emph{%
image} under $R_{\mathbb{T}}$, i.e.,
\begin{equation*}
K_{\mathbb{T}}(A)=\left\{ R_{\mathbb{T}}((d_k))|\mbox{ }A\left(
d_{k},d_{k+1}\right) =1\mbox{ for }k\geq 1\right\} \mbox{.}
\end{equation*}%
We say $K_{\mathbb{T}}(A)$ is a \emph{Sierpi\'{n}ski carpet} and denote by $%
Z=Z_{\left( \mathbf{m},\mathbf{n}\right) }(A)$ the Sierpi\'{n}ski carpet for a
given pair $\left( \mathbf{m},\mathbf{n}\right) $ and transition matrix $A$.
To avoid confusion we will call this a \emph{Markov Sierpi\'{n}ski carpet}.
McMullen \cite{McMullen1984} computes the Hausdorff dimension according to the following
formula.

\begin{theorem}[\cite{McMullen1984}]
\label{Thm: 9}Let $Z$ be a Markov Sierpi\'{n}ski carpet. Construct $\mathbf{n}$
matrices $A_{1},\cdots ,A_{\mathbf{n}}$ which are indexed by $D$ as follows:
$A_{j}(d,d^{\prime })=A(d,d^{\prime })$, if the second coordinate of $%
d^{\prime }\in D$ is $j$ and $A_{j}(d,d^{\prime })=0$, otherwise. Then the
Hausdorff dimension of $Z$ can be formulated as the following formula.
\begin{equation}
\dim _{H}Z=\frac{1}{\log \mathbf{n}}\lim_{k\rightarrow \infty }\frac{1}{k}%
\log \sum_{1\leq i_{0},\cdots ,i_{k-1}\leq \mathbf{n}}\left\Vert
A_{i_{0}}\cdots A_{i_{k-1}}\right\Vert ^{\alpha },  \label{5}
\end{equation}%
where $\alpha =\frac{\log \mathbf{n}}{\log \mathbf{m}}\leq 1$.
\end{theorem}

Notably, Kenyon and Peres \cite{Kenyon1996} extend formula (\ref{5})
to a sofic Sierpi\'{n}ski carpet. In view of (\ref{5}), since
$\left\{A_1, \ldots, A_{\mathbf{n}}\right\}$ are collections of matrices, our
goal in this investigation is to look more closely at how it
relates to the thermodynamic formalism with the matrix-valued
potential function. On the other hand, the more we know about the structure of
$\left\{A_1, \ldots, A_{\mathbf{n}}\right\}$ also enables us to establish more information about the projection space (defined later).

First, we recall the results of Olivier \cite{Olivier2009} in the study of the full
Hausdorff dimension of sofic or Markov Sierpi\'{n}ski carpets. Let $Z$ be a
Markov or sofic shift, let $\sigma _{Z}:Z\rightarrow Z$ be its shift map, and the
author defines the so-called (\textbf{H}) condition on $\sigma _{Z}$. $%
\sigma _{Z}$ is said to satisfy the condition \textbf{(H)} if the $y$-axis projection $%
\pi _{y}\mu $ of the Parry measure $\mu $ on $Z$ is a $\phi $-conformal
measure of some normalized potential $\phi :Y\rightarrow \mathbb{R}$. The
condition therein was imposed to ensure that the invariant measures of the full
Hausdorff dimension are the equilibrium states of some potential function,
and the Hausdorff dimension formula (\ref{17}) on $Z$ holds.%
\begin{equation}
\dim _{H}Z=\frac{h_{top}(Z)}{\log \mathbf{m}}+\frac{P(\sigma _{Y},\alpha
\phi )}{\log \mathbf{n}}\mbox{.}  \label{17}
\end{equation}

In the first part of this paper, we define an \emph{induced
matrix-valued potential }$N$ on $Y$, and present the criterion for the
existence and uniqueness for the full Hausdorff dimension on $Z$. Meanwhile,
we derive the analogous formula for the Hausdorff dimension (see (\ref{10})). We
emphasize here that the (\textbf{H}) condition may not be satisfied in our
assumption, however, (\ref{10}) still holds. That is, (\ref{10}) holds under
extensive conditions, namely, the irreducibility of the induced
matrix-valued potential function.

Before formulating our main first result, we give some definitions. Let $(%
\mathbf{m},\mathbf{n})\in \mathbb{Z}_{+}^{2}$ and define two orders $\prec
_{x}$ and $\prec _{y}$ on $D$: we say $d\prec _{x}d^{\prime }$ if $%
d_{1}<d_{1}^{\prime }$ or $d_{1}=d_{1}^{\prime }$ and $d_{2}<d_{2}^{\prime }$%
. According to this order, every $d\in D$ has a unique number on $\left\{
1,\ldots ,\mathbf{mn}\right\} $, we denote by $\Xi ^{(x)}:D\rightarrow
\left\{ 1,\ldots ,\mathbf{mn}\right\} $ the map which assigns each element in
$D$ to the unique number in $\left\{ 1,\ldots ,\mathbf{mn}\right\} $. Define
the order $\prec _{y}$ in the same fashion: $d\prec _{y}d^{\prime }$ if $%
d_{2}<d_{2}^{\prime }$ or $d_{2}=d_{2}^{\prime }$ and $d_{1}<d_{1}^{\prime }$%
. Let $\Xi ^{(y)}:D\rightarrow \left\{ 1,\ldots ,\mathbf{mn}\right\} $ be
also defined similarly, we denote by $\tau _{\left( \mathbf{m},\mathbf{n}%
\right) }$ the permutation on $D$: $\tau \left( d\right) =d^{\prime }$ if $%
\Xi ^{(x)}\left( d\right) =\Xi ^{(y)}\left( d^{\prime }\right) $ and denote
by $P_{\left( \mathbf{m},\mathbf{n}\right) }$ the permutation matrix
associated with the permutation $\tau _{\left( \mathbf{m},\mathbf{n}\right)
} $.

Let $P=P_{\left( \mathbf{m},\mathbf{n}\right) }$, define $B=PAP^{-1}$ and
regard $B$ as $\mathbf{n}\times \mathbf{n}$ system with entries are $\mathbf{%
m}\times \mathbf{m}$ matrices. That is\footnote{%
We note here that we use the index form of $A^{(i,j)}$ to denote the $\left(
i,j\right) $-coordinate of $A$ and $A^{(i,j)}$ is a matrix. And we use the
standard form $A\left( i,j\right) $ to denote the $\left( i,j\right) $%
-coordinate of $A$ if it is a real value.},
\begin{eqnarray}
B &=&\left[
\begin{array}{ccc}
B^{\left( 1,1\right) } & \cdots & B^{\left( 1,\mathbf{n}\right) } \\
\vdots & \ddots & \vdots \\
B^{\left( \mathbf{n},1\right) } & \cdots & B^{\left( \mathbf{n},\mathbf{n}%
\right) }%
\end{array}%
\right]  \nonumber \\
&=&\left[
\begin{array}{ccc}
\left( PAP^{-1}\right) ^{\left( 1,1\right) } & \cdots & \left(
PAP^{-1}\right) ^{\left( 1,\mathbf{n}\right) } \\
\vdots & \ddots & \vdots \\
\left( PAP^{-1}\right) ^{\left( \mathbf{n},1\right) } & \cdots & \left(
PAP^{-1}\right) ^{\left( \mathbf{n},\mathbf{n}\right) }%
\end{array}%
\right]  \label{27}
\end{eqnarray}%
The matrix-valued potential function $N=\left( N_{ij}\right) _{i,j=1}^{%
\mathbf{n}}$ arising from (\ref{27}) is defined:
\begin{equation}
N=\left( N_{ij}\right) _{i,j=1}^{\mathbf{n}}=\left( B^{\left( i,j\right)
}\right) _{i,j=1}^{\mathbf{n}}=\left( \left( PAP^{-1}\right) ^{\left(
i,j\right) }\right) _{i,j=1}^{\mathbf{n}}.  \label{28}
\end{equation}%
{}

We adapt the name from \cite{Chazottes2003, Chazottes2011} to call $N$ the \emph{%
induced (matrix-valued) potential }on $Y$. The \emph{normalized} \emph{%
induced (matrix-valued) potential }$\bar{N}=\left( \bar{N}_{ij}\right)
_{i,j=1}^{\mathbf{n}}$ is also defined by $\bar{N}_{ij}=\rho _{A}^{-1}N_{ij}$
for all $1\leq i, j \leq \mathbf{n}$, where $\rho _{A}$ denotes the maximal
eigenvalue of matrix $A$.

A family of $n\times n$ matrices $\left( N_{i}\right) _{i\in \mathcal{S}}$
with entries in $\mathbb{R}$ is said to be \emph{irreducible} over $\mathbb{R%
}^{n}$ if there is no non-zero proper linear subspace $V$ of $\mathbb{F}^{n}$
such that $N_{i}V\subseteq V$ for all $i\in \mathcal{S}$. The first result
of this investigation is the following.

\begin{theorem}
\label{Thm: A_Rev: 1}Let $Z$ be a Markov Sierpi\'{n}ski carpet and $N=\left(
N_{ij}\right) _{i,j=1}^{\mathbf{n}}$ be the induced potential from $A$.
Assume $N$ is irreducible, then

\begin{enumerate}
\item[(i)] The following statements are equivalent.

\begin{enumerate}
\item[(a)] $\mu $ is the unique measure of the full Hausdorff dimension.

\item[(b)] $\mu $ is the unique conditional equilibrium measure (defined in Section 3.1) of the zero
potential function on $Z$ with respect to $\nu _{\alpha }$, where
$\nu _{\alpha }$ is the unique equilibrium measure of the
matrix-valued potential $\alpha N=\left( \| N_J\|^{\alpha}
\right)_{J\in Y^*}.$
\end{enumerate}

\item[(ii)] The following Hausdorff dimension formula holds:%
\begin{equation}
\dim _{H}Z=\frac{h_{top}(Z)}{\log \mathbf{m}}+\frac{P(\sigma _{Y},\alpha N)}{%
\log \mathbf{n}}  \label{10}
\end{equation}
\end{enumerate}
\end{theorem}

The essential ingredient of the proof in Theorem \ref{Thm: A_Rev: 1}
is that the irreducibility of $N$ ensures the existence of the Gibbs
measure $\nu $. Since the Gibbs measure $\nu $ may have
\emph{infinite memory} (cf.~\cite{Chazottes2003, Chazottes2011}),
the question arises: which conditions ensure that the measure $\nu $
has \emph{finite memory}? The structure of a \emph{$k$-th higher
block induced (matrix-valued) potential} $N^{[k]}$ (defined later)
plays an important role in answering this question. We denote by $Y_k$ the collection of all possible
words in $Y$ of length $k$. For $k \geq 2,$ let
\begin{equation*}
D^{[k]} = \left\{\left(d_1^{(1)}\ldots d_1^{(k)}, d_2^{(1)}\ldots
d_2^{(k)}\right): (d_1^{(i)}, d_2^{(i)})\in D \hbox{ for all } i= 1,
\ldots, k\right\},
\end{equation*}
$d^{[k]} = (d^{[k]}_1, d^{[k]}_2) = \left(d_1^{(1)}\ldots d_1^{(k)}, d_2^{(1)}\ldots
d_2^{(k)}\right)$ and $A^{[k]}\in \mathbb{R}^{\textbf{d}^k \times \textbf{d}^k}
(\textbf{d} = \textbf{m} \times \textbf{n})$ be the $k$-th higher
block transition matrix from $A$ which is indexed by $D^{[k]}$. We
define the permutation matrix $P^{[k]}= P^{[k]}_{(\textbf{m},
\textbf{n})}$ in the same fashion as $P = P^{[1]}_{(\textbf{m},
\textbf{n})}$. Set $B^{[k]} =
P^{[k]}A^{[k]}\left(P^{[k]}\right)^{-1}$ and regards $B^{[k]}$ as
$\textbf{n}^k \times \textbf{n}^k$ system with entries are
$\textbf{m}^k \times \textbf{m}^k$ matrices. The matrix-valued
potential function $N^{[k]} = \left(N^{[k]}_J\right)_{J\in Y_{k+1}}$
is defined by
\[
N^{[k]}_J  = \left( B^{[k]}\right)^{\left(J(0, k-1),
J(1, k)\right)} \hbox{ for all } J= J(0,k)\in Y_{k+1}.
\]
We call $N^{[k]}$ the \emph{k-th higher block induced
(matrix-valued) potential} on $Y^{[k]}.$ Note that $N=N^{[1]}$, and $N^{[i]}$ is defined by $A^{[i]}$ and $P^{[i]}$, for $i=1, \ldots, k.$ (Figure \ref{fig1})

\begin{center}
\label{fig1}
\begin{pspicture}(8,3)
\rput(1,2.5){\rnode{a}{$A^{[1]} =:$}}\rput(2,2.5){$A$}
\rput(1,0.5){\rnode{A}{$N^{[1]} =:$}}\rput(2,0.5){$N$}
\rput(2,2.2){\rnode{A}{}}\rput(2,0.8){\rnode{N}{}}\ncline{->}{A}{N}\Bput{$P$}
\rput(4,2.5){$A^{[2]}$}\rput(6,2.5){$\cdots$}\rput(8,2.5){$A^{[k]}$}
\rput(4,0.5){$N^{[2]}$}\rput(8,0.5){$N^{[k]}$}
\rput(4,2.2){\rnode{B}{}}\rput(4,0.8){\rnode{C}{}}\ncline{->}{B}{C}\Bput{$P^{[2]}$}
\rput(8,2.2){\rnode{D}{}}\rput(8,0.8){\rnode{E}{}}\ncline{->}{D}{E}\Bput{$P^{[k]}$}
\rput(2.3,2.5){\rnode{f}{}}\rput(3.6,2.5){\rnode{g}{}}\ncline{->}{f}{g}
\rput(4.4,2.5){\rnode{i}{}}\rput(5.6,2.5){\rnode{j}{}}\ncline{->}{i}{j}
\rput(6.4,2.5){\rnode{k}{}}\rput(7.6,2.5){\rnode{l}{}}\ncline{->}{k}{l}
\end{pspicture}
\figcaption{For every matrix $A^{[i]}$, there exists a permutation matrix $P^{[i]}$ such that the induced matrix $N^{[i]}$ is obtained by applying $P^{[i]}$ on $A^{[i]}$.}
\end{center}

If $J\in Y_{k}$ and $0\leq m, n\leq k$, we use the notation $%
J\left( m,n\right) $ to denote the subword of $J$ from coordinate $m$ to $n$%
, i.e., $J(m,n)=\left( j_{m}\ldots j_{n}\right) $ if $J=(j_{0}j_{1}%
\ldots j_{k})$. For $k\geq 1$ and $J\in Y_n$ with $n \geq k$, $N^{[k]}_J$ stands for the product of matrices of $N^{[k]}$ along the
path of $J$, i.e.,
\[
N^{[k]}_J = N^{[k]}_{J(1,n)} = \prod_{i=1}^{n-k}N^{[k]}_{J(i, i+k)}
\hbox{ for all } J=J(1,n)\in Y_n.
\]

We say that $N$ satisfies the \emph{Markov condition from
left of order }$k$ if there exists a non-zero linear subspace
$\left\{V_J\right\}_{J\in Y_k}\subseteq \mathbb{R}^{\textbf{m}^k}$ such that $%
V_{J\left( 0,k-1\right) }N^{[k]}_{J\left( 0,k\right) }\subseteq
V_{J\left( 1,k\right) }$ for all $J\left( 0,k-1\right) $ and
$J\left( 1,k\right) \in Y_{k}$, $N$ satisfies the \emph{Markov
condition from right of order k} if $\left(N^{[k]}\right)^{t}=\left(
(N^{[k]}_J)^{t}\right) _{J\in Y_{k+1}}$ satisfies the Markov
condition from left, where $A^{t}$ denotes the transpose of $A$.
Finally, say $N$ satisfies the \emph{Markov condition} if it
satisfies either the Markov condition from the left or right for
some order $k\in \mathbb{N}$.
The following theorem provides a criterion for
checking whether $\nu $ is a Markov measure.

\begin{theorem}
\label{Thm: 8}Let $Z$ be a Markov Sierpi\'{n}ski carpet and $N =
\left(N_{ij}\right)_{i, j=1}^{\mathbf{n}}$ be the induced
matrix-valued potential from $A$. Then, $\nu $ is a $k$-step Markov
measure on $Y$ if and only if $N$
satisfies the Markov condition of order $k$. Furthermore, if $\nu $ is a $k$%
-step Markov measure, then $k\leq \mathbf{m}-\mathbf{n}$.
\end{theorem}

We mention here that the inequality $k\leq \mathbf{m}-\mathbf{n}$ in Theorem %
\ref{Thm: 8} is sharp.  More precisely, we examine the well-known example of a McMullen carpet (Example \ref{Ex: 1}) in which the induced matrix-valued
potential
satisfies the Markov condition of order $1$. It also follows from Theorem %
\ref{Thm: 8} and the fact that $\mathbf{m}=3$ and $\mathbf{n}=2$, that the
Markov measure induced from $N$ can only be $1$-step.

Compared to Theorem \ref{Thm: A_Rev: 1}, Theorem \ref{Thm: 8} reveals
that the more structured the vector space of $V_{J}$ from
$N^{[k]}$, the more it implies about the property of $\nu $. In other words,
Theorem \ref{Thm: 8} illustrates that the projection measure $\nu $
is Markov if and only if the collection of $\left(
N^{[k]}V_{J}\right) _{J\in Y_{k}}$ is a finite set, which guarantees that
the Gibbs measure $\nu $ falls into the finite range. (Readers may consult \cite{Chazottes2003, Chazottes2011} for more detail.)

At this point, a further question arises: \emph{If }$N$\emph{\
satisfies the Markov condition of order k}\emph{, what kind of Markov measure
is }$\nu $\emph{?} To answer this question, we may assume $m\left(
J\left( 0,k-1\right) ,J\left(
1,k\right) \right) \in \mathbb{R}$ such that%
\begin{equation*}
V_{J\left( 0,k-1\right) }N^{[k]}_{J\left( 0,k\right) }=m\left(
J\left( 0,k-1\right) ,J\left( 1,k\right) \right) V_{J\left(
1,k\right) },
\end{equation*}%
the following theorem illustrates that the coefficient of $m\left(
J,J^{\prime }\right) $ helps us to determine what kind of Markov measure $%
\nu $ is.

\begin{theorem}
\label{Thm: 3}If $N$ satisfies the Markov condition of order $k$, then $\nu
$ is the unique maximal measure of the subshift of finite type $X_M$
with adjacency matrix $M=\left[ m\left( J,J^{\prime }\right) \right]
_{J,J^{\prime }\in Y_{k}}$.
\end{theorem}

Let us return to the Markov or sofic Sierpi\'{n}ski carpet. We
recall the two following interesting problems:

\begin{itemize}
\item[(i)] When are the Hausdorff and Minkowski dimensions coincident?

\item[(ii)] What is the exact value of the Hausdorff dimension?
\end{itemize}

These two problems seems to have satisfactory answers when $Z$ is a Markov
Sierpi\'{n}ski carpet. For (i), Kenyon and Peres show that if $A$ is primitive,
then $\dim _{H}Z=\dim _{M}Z$ if and only if the unique invariant measure of
maximal entropy on $Z$ projects via $\pi _{y}$ to the unique measure of
maximal entropy on $\pi _{y}(Z)$ (Theorem \ref{Thm: 13}). For problem (ii), if $Z$ is a Markov
Sierpi\'{n}ski carpet, let $D^{\prime }\subseteq D=\left\{ 1,\ldots ,\mathbf{m}%
\right\} \times \left\{ 1,\ldots ,\mathbf{n}\right\} $ be the non-empty subset
of $D$. Define
\begin{equation}
K(T,D^{\prime })=\left\{ \sum_{k=1}^{\infty }
\left(\begin{array}{cc}
\mathbf{m}^{-k} &0 \\
0 &\mathbf{n}^{-k}
\end{array} \right)
d_{k}:d_{k}\in
D^{\prime }\mbox{ for all }k\right\} \mbox{.}  \label{12}
\end{equation}%
The Hausdorff dimension of $Z^{\prime }=K(T,D^{\prime })$ has a closed
form: let $z(j)$ be the number of rectangles in row $j$
\begin{equation}
\dim _{H}Z^{\prime }=\frac{1}{\log \mathbf{n}}\log \sum_{j=1}^{\mathbf{n}%
}z(j)^{\alpha },\mbox{ where }\alpha =\frac{\log \mathbf{n}}{\log \mathbf{m}}\mbox{.}
\label{15}
\end{equation}

In the following, the structure of $N$ helps us to derive the closed formula
for a more general Sierpi\'{n}ski carpet and the explicit value for the Hausdorff
dimension. Assume that $N$ satisfies the Markov condition of order $1$. Define the \emph{induced graph} and the corresponding \emph{induced} \emph{%
transition matrix} as follows: Let $T$ be given and $N$ be the induced
matrix-valued potential from $T$, then let $\mathcal{V} = \left\{ 1,\ldots ,%
\mathbf{n}\right\} $ and $\mathcal{E}=\left\{ (i,j)\right\} _{i,j=1}^{%
\mathbf{n}}$, where $(i,j)=1$ if $N_{ij}$ is non-zero matrix. We call $%
G=\left( \mathcal{V},\mathcal{E}\right) $  the induced graph. Define
\begin{equation*}
T_G(i, j) =\left\{
             \begin{array}{ll}
               1, & \hbox{if $N_{ij} \neq 0_{\mathbf{m}\times \mathbf{m}}$;} \\
               0, & \hbox{otherwise}
             \end{array}
           \right.
\end{equation*}
the induced transition matrix corresponds to $G$. Finally, we define
$G^{[k]}$ and $T^{[k]} := T_{G^{[k]}}$ in the same fashion if $N$
satisfies the Markov condition of order $k >1$.

\begin{theorem}
\label{Thm: 6}Let $Z$ be a Markov Sierpi\'{n}ski carpet. Assume $N$ satisfies the Markov condition of order
$k$ and let $T^{[k]}$ be the induced transition matrix which corresponds to induced graph $G^{[k]}$. Then

\begin{enumerate}
\item[(i)] $\rho_{M}=\rho_{T^{[k]}}$ if and only if $\dim _{H}Z=\dim
_{M}Z $, where $M$ is defined in Theorem \ref{Thm: 3}.

\item[(ii)] Let $D^{\prime }\subseteq D$ and $K(T,D^{\prime })$ as defined in
(\ref{12}), then (\ref{15}) holds.

\item[(iii)] Define $M^{\alpha }=\left[ m^{\alpha}\left( J, J^{\prime }\right) \right] _{J, J^{\prime }\in
Y_{k}}$. Then $\dim _{H}Z = \log_{\mathbf{n}} \rho _{M^{\alpha }},$ where $\rho
_{M}$ is the maximal eigenvalue of $M$ and $\alpha = \log\textbf{n}/\log\textbf{m}$.
\end{enumerate}
\end{theorem}

\begin{remark}
If $N=\left( N_{ij}\right) _{i,j=1}^{\mathbf{n}}$ is reducible, Proposition
1.4 of \cite{FengKaenmaki2011} demonstrates that one can decompose $N$ to the
irreducible components. This reveals that the equilibrium measures for $N$
may not be unique. On the other hand, it can be easily checked whether or not the
reducibility of $N=\left( N_{ij}\right) _{i,j=1}^{\mathbf{n}}$ implies the
reducibility of $A$ (Since $N$ is extracted from $B$ which is the
permutation of $A$). This illustrates that the non-uniqueness for the
equilibrium measure of $\nu _{\alpha }$ on $Y$ relates to the non-uniqueness
for the maximal measure on $Z$ of $A$.
\end{remark}

The rest of the paper is organized as follows. Since the space $Y$ with the
induced potential $N$ is no longer $p$-specification (cf.~\cite{BarralFeng2009, Falconer1988, Falconer2003}), it is weak $p$-specification instead. We review some known results in \cite{Feng2011} for weak $p$-specification shift in Section 2. The detailed proofs for Theorem \ref{Thm: A_Rev: 1}, Theorem \ref{Thm: 8} and Theorem \ref{Thm: 3} are
presented in Section 3.

In Section 4, the established results for the induced potential $N$ enable
us to answer problems raised by Chazottes and Ugaldes \cite{Chazottes2003}, and Boyle
and Petersen \cite{Boyle2011}. To be precise, Chazottes and Ugaldes use the
ansatz of the induced potential to prove the existence of well-defined
potential function, and the corresponding Gibbs measure (BGM \cite{Chazottes2003})
on the projection space under (\textbf{H1}) and (\textbf{H2}). They raise
the following problem: \emph{When is the factor map not a topological Markov
map?} On the other hand, Boyle and Petersen raise the following question ({\cite[Problem
3.3]{Boyle2011}}): \emph{Given a procedure to decide, and given a factor map }$%
\pi :\Omega _{A}\rightarrow \Omega _{B}$, \emph{where $\Omega_{A}$ and
$\Omega_B$ are the Markov system induced by the transition matrices
$A$ and $B$}\emph{, how can we know whether }$\pi $\emph{\ is Markovian?} Theorem
\ref{Thm: 12} is presented in Section 4 to provide a criterion for determining whether such $\pi $ is Markovian. Finally, we also list some interesting
examples, namely, the Blackwell and McMullen examples therein.

\section{Preliminaries}

Let $Z$ be a Markov Sierpi\'{n}ski carpet introduced in Section 1, and define the
sliding block code $\Pi _{y}:D\rightarrow \left\{ 0,\ldots ,\mathbf{n}%
-1\right\} $ by
\begin{equation*}
\Pi _{y}(d)=d_{2}\mbox{ if }d=\left( d_{1},d_{2}\right) \in D\mbox{.}
\end{equation*}%
Denote by $G_{A}$ the graph associated with the adjacent matrix $A$. Then
the pair $\mathcal{G}=\left( G_{A},\Pi _{y}\right) $ forms a one-block
factor map from alphabet $Z$ to $Y=\left( \Pi _{y}\right) _{\infty }(Z)$ as
follows:%
\begin{equation}
\left( \Pi _{y}\right) _{\infty }(Z)=\left\{ y\in \left\{ 0,\ldots ,\mathbf{n%
}-1\right\} ^{\mathbb{N}}:\mbox{ }y=\left( \Pi _{y}\right) _{\infty }(z)%
\mbox{ for some }z\in Z\right\} \mbox{,}  \label{13}
\end{equation}%
where $\left( \Pi _{y}\right) _{\infty }(z)=\Pi _{y}(z_{0})\Pi
_{y}(z_{1})\ldots \in \left\{ 1,\ldots ,\mathbf{n}\right\} ^{\mathbb{N}}$.
In the following, we write $\pi _{y}$ instead of $\left( \Pi _{y}\right)
_{\infty }$.

We say that $X$ satisfies the criterion for \emph{weak specification }\cite{Feng2011} if there
exists $p\in \mathbb{N}$ such that, for any two words $I$ and $J\in X^{\ast
}=\cup _{n\in \mathbb{N}}X_{n}$, where $X_{n}$ is collection admissible
words in $X$ of length $n$, there is a word $K$ of length not exceeding $p$
such that the word $IKJ\in X^{\ast }$.

Denote by $\mathcal{D}_{w}(X,p)$ \cite{Feng2011} the collection of functions $%
f:X^{\ast }\rightarrow \lbrack 0,\infty )$ such that $f\left( I\right) >0$
for at least one $I\in X^{\ast }$ and there exists $0<c\leq 1$ so that

\begin{enumerate}
\item[(1)] $f\left( IJ\right) \leq c^{-1}f\left( I\right) f\left( J\right) $ for all $I, J \in X^{\ast }.$

\item[(2)] For all $I$, $J\in X^{\ast },$ there exists $K\in \cup
_{i=0}^{p}X_{i}$ such that $IKJ\in X^{\ast }$ and $f\left( IKJ\right) \geq
cf\left( I\right) f\left( J\right) .$
\end{enumerate}

One can easily check that $\mathcal{D}_{w}(X,p)\neq \emptyset $ if
and only if $X$ satisfies the weak $p$-specification. Let $N=\left(
N_{ij}\right) _{i,j=1}^{\mathbf{n}}$ be the induced potential from
$A$ on $Y$ and $f(J)=\left\Vert N_{J}\right\Vert $, by Proposition
2.8 and Lemma 2.1 of \cite{Feng2009}, we see that $f\in
\mathcal{D}_{w}(Y,p)$. It also follows from Theorem 5.5
\cite{Feng2011}, $\alpha N = \left(f(J)^{\alpha }\right)_{J\in Y^*}$
has a unique equilibrium $\nu _{\alpha }$. Finally, Theorem 6.1 of \cite%
{Feng2011} shows that if $\nu _{\alpha }$ is the unique equilibrium measure
of $\alpha N$\textbf{,} the zero potential function on $Z$ has a unique $%
\mathbf{a}$-weighted equilibrium state which is the conditional equilibrium
states of $\nu _{\alpha }$ with respect to $\Phi $. We present some useful
Lemmas as follows.

\begin{lemma}[{\cite[Lemma 5.2]{Feng2011}}]
\label{Lma: 2}Suppose $f\in \mathcal{D}_{w}(X,p)$. Then the following two
properties hold:

\begin{enumerate}
\item[(i)] There exists a constant $\gamma >0$ such that for each $I\in
X^{\ast }$, there exists $i,j\in \mathcal{A}(X)$ such that $f(iI)\geq \gamma
f\left( I\right) $ and $f\left( Ij\right) \geq \gamma f\left( I\right) $,
where $\mathcal{A}(X)$ denotes the symbol set on $X$.

\item[(ii)] Let $u_n=\sum_{J\in X_{n}}f(J).$ Then the limit $%
u=\lim_{n\rightarrow \infty }\left( 1/n\right) \log u_{n}$ exists and $%
u_{n}\approx \exp \left( nu\right) $.
\end{enumerate}
\end{lemma}

\begin{lemma}
\label{Lma: 1}Let $Z=Z_{\left( \mathbf{m},\mathbf{n}\right) }(A)$ with $A\in
\mathbb{R}^{\mathbf{d}\times \mathbf{d}}$ be irreducible, where $\mathbf{d}=%
\mathbf{m}\times \mathbf{n}$. Let $A_{1},\ldots ,A_{\mathbf{n}}$ be as
defined in Theorem \ref{Thm: 9} and $B_{k}=PA_{k}P^{-1}$, for $k=1,\ldots ,%
\mathbf{n}$, and we write $B_{k}=\left( B_{k}^{\left( i,j\right) }\right)
_{i,j=1}^{\mathbf{n}}$. Then $B_{k}^{(i,j)}=N_{ij}$, if $j=k$ and $1\leq
i\leq \mathbf{n}$\textbf{,} and $B_{k}^{(i,j)}=0_{\mathbf{m}\times \mathbf{m}%
}$ otherwise, where $0_{\mathbf{m}\times \mathbf{m}}$ denotes the $\mathbf{m}%
\times \mathbf{m}$ matrix with all entries being $0^{\prime }$s.
\end{lemma}

\begin{proof}
It follows from the definition of $N$ and the permutation $P$, it can be
easily checked that the index of the matrix $B=\left( B^{\left( i,j\right)
}\right) _{i,j=1}^{\mathbf{n}}$ equals to the $D$ with the order $\prec _{y}$%
. Define
\begin{equation*}
E_{d_{2}}=\left\{ d_{1}:\left( d_{1},d_{2}\right) \in D\mbox{ with }\Pi
_{y}(d)=d_{2}\right\} \mbox{.}
\end{equation*}%
We see that $B^{\left( i,j\right) }$ is indexed by $E_{i}\times E_{j}$ and $%
B^{\left( i,j\right) }\left( p, q\right) =1$ if $A\left( \left( p,i\right)
,\left( q, j\right) \right) =1$ with $p \in E_{i}$ and $q \in E_{j}$. Take $%
1\leq k\leq \mathbf{n}$, it follows from the definition of $A_{k}$: $%
A_{k}(d,d^{\prime })=A(d,d^{\prime })$ if the second coordinate of $%
d^{\prime }\in D$ is $k$, it means that $B_{k}^{(i,j)}$ is indexed by $%
E_{i}\times E_{j}$ for which $B_{k}^{(i,j)}\left( p, q\right) =1$ if and only
if $A\left( (p, i),(q, j)\right) =1$ and $j= k.$ Therefore,
\begin{equation*}
B_{k}^{(i,j)}\left( p, q\right) =\left( PA_{k}P^{-1}\right) \left( p, q\right)
=N_{ik}\mbox{ for all }1\leq i\leq \mathbf{n}\mbox{\textbf{.}}
\end{equation*}%
The proof is thus completed.
\end{proof}

\begin{theorem}
\label{Thm: 1}Let $Z=Z_{(\mathbf{m},\mathbf{n})}\left( A\right) $ be a Markov
Sierpi\'{n}ski carpet with $A$, assume that $N$ the induced potential from $A$ is irreducible. Then,%
\begin{equation}
\dim _{H}Z=\frac{1}{\log \mathbf{n}}\lim_{n\rightarrow \infty }\frac{1}{n}%
\log \sum_{J\in Y_{n}}\left\Vert N_{J}\right\Vert ^{\alpha }  \label{6}
\end{equation}%
where $\alpha =\log \mathbf{n}/\log \mathbf{m}$.
\end{theorem}

\begin{proof}
Let $Z=Z_{\left( \mathbf{m},\mathbf{n}\right) }(A)$ be given and $\pi
_{y}:Z\rightarrow Y$ be a sliding block code from $Z$ to $Y$ as in (\ref{13}%
) and recall that $\mathbf{d}=\mathbf{m}\times \mathbf{n}$. We first show
that there exists $c>0$ such that for all $k\in \mathbb{N}$ and $%
J=(j_{0},\ldots ,j_{n-1})\in Y_{n}$
\begin{equation}
c^{-1}\left\Vert N_{j_{0}j_{1}}\ldots N_{j_{n-2}j_{n-1}}\right\Vert \leq
\left\Vert B_{J}\right\Vert =\left\Vert A_{J}\right\Vert \leq c\left\Vert
N_{j_{0}j_{1}}\ldots N_{j_{n-2}j_{n-1}}\right\Vert ,  \label{11}
\end{equation}%
where $B_{J}=B_{j_{0}}B_{j_{1}}\cdots B_{j_{n-1}}$ and $B_{k}$ is
defined in Lemma \ref{Lma: 1}. Indeed, since $P$ is a permutation,
we have
\begin{eqnarray}
\left\Vert A_{J}\right\Vert
&=&\mathbf{1}_{\textbf{d}}^{t}A_{j_{0}}A_{j_{1}}\cdots
A_{j_{n-1}}\mathbf{1}_{\textbf{d}}=\mathbf{1}_{\textbf{d}}^{t}P^{-1} B_{j_{0}}B_{j_{1}}\cdots
B_{j_{n-1}}P\mathbf{1}_{\textbf{d}}  \nonumber \\
&=&\mathbf{1}_{\textbf{d}}^{t}B_{j_{0}}B_{j_{1}}\cdots B_{j_{n-1}}\mathbf{1}%
_{\textbf{d}}=\left\Vert B_{J}\right\Vert .  \label{34}
\end{eqnarray}%
Therefore
\begin{equation}
\left\Vert A_{J}\right\Vert =\left\Vert B_{J}\right\Vert \mbox{ for all }%
J\in Y_{n}.  \label{19}
\end{equation}%
Since $A$ is irreducible we conclude that $B$ is also irreducible, we have $%
c_{1}=\max\limits_{0\leq u,v\leq \mathbf{m}-1}\max\limits_{_{0\leq l\leq
\mathbf{n}-1}}\left( \sum_{i=1}^{\mathbf{n}}N_{il}\right) (u,v)>0$.
According to Lemma \ref{Lma: 1} we also have
\begin{eqnarray}
\left\Vert B_{J}\right\Vert
&=&\mathbf{1}_{\textbf{d}}^{t}B_{j_{0}}B_{j_{1}}\cdots
B_{j_{n-1}}\mathbf{1}_{\textbf{d}}=\mathbf{1}_{\textbf{d}}^{t}\left[
\begin{tabular}{lllll}
$0$ & $\cdots $ & $N_{1j_{0}}$ & $\cdots $ & $0$ \\
$\vdots $ & $\ddots $ & $N_{2j_{0}}$ & $\cdots $ & $0$ \\
$\vdots $ & $\cdots $ & $\vdots $ & $\cdots $ & $0$ \\
$0$ & $\cdots $ & $N_{\mathbf{n}j_{0}}$ & $\cdots $ & $0$%
\end{tabular}%
\right] \times \cdots  \nonumber \\
&&\times B_{j_{n-2}}\times \left[
\begin{tabular}{lllll}
$0$ & $\cdots $ & $N_{1j_{n-1}}$ & $\cdots $ & $0$ \\
$\vdots $ & $\ddots $ & $N_{2j_{n-1}}$ & $\cdots $ & $0$ \\
$\vdots $ & $\cdots $ & $\vdots $ & $\cdots $ & $0$ \\
$0$ & $\cdots $ & $N_{\mathbf{n}j_{n-1}}$ & $\cdots $ & $0$%
\end{tabular}%
\right] \mathbf{1}_{\textbf{d}}  \nonumber \\
&=&\mathbf{1}_{\mathbf{m}}^{t}\left( \sum_{i=1}^{\mathbf{n}%
}N_{ij_{0}}\right) N_{j_{0}j_{1}}\cdots N_{j_{n-2}j_{n-1}}\mathbf{1}_{%
\mathbf{m}}^{t}  \label{35}
\end{eqnarray}%
Therefore,
\begin{equation}
\left\Vert B_{J}\right\Vert =\mathbf{1}_{\textbf{m}}^{t}\left( \sum_{i=1}^{\mathbf{n}%
}N_{ij_{0}}\right) N_{j_{0}j_{1}}\cdots N_{j_{n-2}j_{n-1}}\mathbf{1}_{%
\mathbf{m}}\leq c_{1}\left\Vert N_{J}\right\Vert \mbox{.}  \label{18}
\end{equation}%
On the other hand, since $N$ is irreducible, then Lemma \ref{Lma: 2}-(i) is
applied to show that there exists a $1\leq i\leq \mathbf{n}$ such that
\begin{equation}
\left\Vert B_{J}\right\Vert \geq \gamma \mathbf{1}_{\mathbf{m}%
}^{t}N_{j_{0}j_{1}}\cdots N_{j_{n-2}j_{n-1}}\mathbf{1}_{\mathbf{m}}=\gamma
\left\Vert N_{J}\right\Vert .  \label{14}
\end{equation}%
Combining (\ref{14}), (\ref{19}), (\ref{18}), Lemma \ref{Lma: 2}-(ii),
Theorem \ref{Thm: 9} and (\ref{11}) yields
\begin{eqnarray*}
\dim _{H}Z &=&\frac{1}{\log \mathbf{n}}\lim_{n\rightarrow \infty }\frac{1}{n}%
\log \sum_{J\in Y_{n}}\left\Vert A_{J}\right\Vert ^{\alpha } \\
&=&\frac{1}{\log \mathbf{n}}\lim_{n\rightarrow \infty }\frac{1}{n}\log
\sum_{J\in Y_{n}}\left\Vert B_{J}\right\Vert ^{\alpha } \\
&=&\frac{1}{\log \mathbf{n}}\lim_{n\rightarrow \infty }\frac{1}{n}\log
\sum_{J\in Y_{n}}\left\Vert N_{J}\right\Vert ^{\alpha }\mbox{ (Lemma \ref%
{Lma: 2}-(ii)).}
\end{eqnarray*}%
The theorem is thus proved.
\end{proof}

\begin{remark}

\begin{enumerate}
\item We define the limit in (\ref{6}) as the topological pressure $P(\sigma
_{Y},\alpha N)$ on $Y$ with respect to the potential function $N$.

\item We note here that Yayama \cite{Yayama2011} derived a similar result as in
Theorem \ref{Thm: 1} ({\cite[Theorem 4.4-(1)]{Yayama2011}}). To be precise, the
author proves the following:%
\begin{eqnarray}
\dim _{H}Z &=&\frac{1}{\log \mathbf{n}}\lim_{n\rightarrow \infty }\frac{1}{n}%
\log \sum_{J\in Y_{n}}\left\vert \pi _{y}^{-1}\left( y_{1}\cdots
y_{n}\right) \right\vert ^{\alpha }  \label{20} \\
&=&\frac{1}{\log \mathbf{n}}\lim_{n\rightarrow \infty }\frac{1}{n}\log
\sum_{J\in Y_{n}}\left\Vert N_{J}\right\Vert ^{\alpha }\mbox{.}  \nonumber
\end{eqnarray}%
The second equality comes from Lemma \ref{Lma: 1}. We emphasize here that
the potential $\alpha N=f^{\alpha }(J)$ on $Y$ is not necessary
continuous.
\end{enumerate}
\end{remark}

\subsection{Existence and uniqueness of Gibbs measures for matrix-valued
potential}

This section presents the existence and uniqueness of the Gibbs measure for the
matrix-valued potential function with some irreducible condition. Feng and
Kaenmaki \cite{FengKaenmaki2011} characterize the structure of equilibrium and the Gibbs
measure for matrix-valued potentials for irreducible $N.$

\begin{theorem}
\label{Thm: 4}Let $N=\left( N_{i}\right) _{i\in \mathcal{S}}$ be a family of
$d\times d$ matrices with entries in $\mathbb{R}$. If $N=\left( N_{i}\right)
_{i\in \mathcal{S}}$ is irreducible. Then for each $\alpha >0$, $P(\sigma
_{Y},\alpha N)$ has a unique $\alpha $-equilibrium measure $\mu _{\alpha }$
which satisfies the Gibbs property: $\forall n\in \mathbb{N}$ and $J\in
Y_{n} $, there exists $c>0$ such that
\begin{equation}
c^{-1}\exp (-nP(\sigma _{Y},\alpha N))\left\Vert N_{J}\right\Vert ^{\alpha}\leq
\mu _{\alpha }\left( \left[ J\right] \right) \leq c\exp (-nP(\sigma
_{Y},\alpha N))\left\Vert N_{J}\right\Vert ^{\alpha}.\mbox{ }  \label{9}
\end{equation}
\end{theorem}

The following theorem illustrates the existence of the Gibbs measure
on $Y$ with respect to the induced potential function $N$. Suppose
$X$ is a shift space, denote by $\mathcal{M}(X, \sigma_{X})$ the
collection of  all $\sigma_{X}$-invariant measures on $X$.

\begin{corollary}
\label{Cor: 1}Under the same assumptions of Theorem \ref{Thm: A_Rev: 1},
then, for any $\alpha >0$, there exists a unique $\alpha $-equilibrium $\nu
_{\alpha }\in \mathcal{M}(Y,\sigma _{Y})$ which satisfies (\ref{9}).
Furthermore, if $\bar{\nu}_{\alpha }\in \mathcal{M}(Y,\sigma _{Y})$ is the equilibrium
measure of $\alpha \bar{N},$ then $\nu _{\alpha }=\bar{\nu}_{\alpha }$.
\end{corollary}

\begin{proof}
Since $N$ is irreducible, the existence of the unique $\alpha $-equilibrium $%
\nu _{\alpha }$ which satisfies (\ref{9}) is the immediate consequence of
Theorem \ref{Thm: 4}. From the definition of $N$ and $\bar{N}$ we know that $%
N$ is irreducible if and only if $\bar{N}$ is irreducible. Then there exists
a unique $\alpha $-equilibrium measure $\bar{\nu}_{\alpha }$ which satisfies
(\ref{9}) for $\alpha \bar{N}$, i.e., $\forall J\in Y_{n}$ and $n\in \mathbb{%
N},$ we have
\begin{equation*}
d^{-1}\exp (-nP(\sigma_{Y},\alpha \bar{N}))\left\Vert \bar{N}_{J}\right\Vert ^{\alpha }\leq \bar{%
\nu}_{\alpha }(\left[ J\right] )\leq d\exp (-nP(\sigma_{Y},\alpha
\bar{N}))\left\Vert \bar{N}_{J}\right\Vert ^{\alpha }\mbox{,}
\end{equation*}%
for some $d>0$. We claim that $\nu _{\alpha }=\bar{\nu}_{\alpha }$. Indeed,
for each $n\in \mathbb{N}$ and $J\in Y_{n}$ we have
\begin{eqnarray}
\nu _{\alpha }(\left[ J\right] ) &\leq &c\exp (-nP(\sigma_{Y},\alpha
N))\left\Vert
N_{J}\right\Vert ^{\alpha }\mbox{ (From (\ref{9}))}  \nonumber \\
&=&c\exp (-nP(\sigma_{Y},\alpha N))\rho _{A}^{-\alpha
n}\left\Vert \bar{N}_{J}\right\Vert
^{\alpha }\mbox{ (since }\bar{N}=\rho _{A}^{-1}N\mbox{)}  \nonumber \\
&=&c\exp (-nP(\sigma_{Y},\alpha \bar{N}))\rho _{A}^{-\alpha n}\rho
_{A}^{\alpha
n}\left\Vert \bar{N}_{J}\right\Vert ^{\alpha }  \nonumber \\
&=&c\exp (-nP(\sigma_{Y},\alpha \bar{N}))\left\Vert \bar{N}_{J}\right\Vert ^{\alpha }  \nonumber \\
&\leq &cd\bar{\nu}_{\alpha }(\left[ J\right] ).  \label{39}
\end{eqnarray}%
Similarly, we have%
\begin{equation}
\nu _{\alpha }(\left[ J\right] )\geq c^{-1}d^{-1}\bar{\nu}_{\alpha }(\left[ J%
\right] )\mbox{.}  \label{40}
\end{equation}%
The claim follows by combining (\ref{39}), (\ref{40}) and the fact that $\nu
_{\alpha }$ and $\bar{\nu}_{\alpha }$ are both ergodic. This completes the
proof.
\end{proof}

\subsection{Sofic measures and linear representable measures}

Let $(X,\sigma _{X})$ and $(Y,\sigma _{Y})$ be subshifts and $\pi
:X\rightarrow Y$ be a sliding block code, each measure $\mu \in \mathcal{M}%
(X,\sigma _{X})$ determines a measure $\pi \mu \in \mathcal{M}(Y,\sigma
_{Y}) $ by
\begin{equation*}
(\pi \mu )(E)=\mu (\pi ^{-1}(E)),\mbox{ }\forall E\subset Y.
\end{equation*}%
If $\mu $ is a Markov measure, then $\pi \mu $ is called a \emph{sofic
measure} (cf.~\cite{Boyle2011}). Let $B\in \mathbb{R}^{d}$ be an irreducible matrix
with spectral radius $\rho _{B}$ and positive right eigenvector $r$, the
\emph{stochasticization} of $B$ is the stochastic matrix
\begin{equation*}
\mathbb{B}:= stoch(B)=\frac{1}{\rho_B }D^{-1}BD,
\end{equation*}%
where $D$ is the diagonal matrix with diagonal entries $D(i,i)=r(i)$. A
measure $\mu $ on $X$ is called \emph{linear representable }with dimension $%
d $ if there exists a triple $(x,P,y)$ with $x\in \mathbb{R}_{+}^{n}$ being a $%
1\times d$ row vector, $y\in \mathbb{R}_{+}^{n}$ is a $d\times 1$ column
vector and $P=\left( P_{i}\right) _{i\in \mathcal{A}(X)}$, where $P_{i}\in
\mathbb{R}^{d\times d}$ such that for all $I=(i_{0},\ldots ,i_{n-1})\in
X_{n} $, the measure $\mu $ can be characterized as the following form:
\begin{equation*}
\mu (\left[ I\right] )=xP_{I}y,
\end{equation*}%
where $P_{I}=P_{i_{0}}P_{i_{1}}\cdots P_{i_{k-1}}$(readers may refer to \cite%
{Boyle2011} for more detail). The triple $(x,P,y)$ is called the
\emph{linear representation} of the measure $\mu$.

\begin{proposition}[{\cite[Theorem 4.20]{Boyle2011}}]
\label{Prop: 3}Let $X=X_{A}$ be a Markov shift with adjacent matrix $A\in
\mathbb{R}^{n\times n}$ which is irreducible and $\pi :X\rightarrow Y$ be a
factor induced from one block map $\Pi :\mathcal{A}(X)\rightarrow \mathcal{A}%
(Y)$, i.e., $\pi =\Pi _{\infty }$. Let $\mathbb{A}=stoch(A)$ and $l$ be the
left eigenvector of $\mathbb{A}$. Then

\begin{enumerate}
\item[(i)] The Markov measure $\mu _{A}$ on $X$ is the linear representable
measure with respect to the triple $\left( x,P,y\right) $, where $x=l,$ $y=%
\mathbf{1}_{n}$, where $P$ is generated by $\left( P_{i}\right) _{i\in
\mathcal{A}(X)}=\left( \mathbb{A}_{i}\right) _{i\in \mathcal{A}(X)}$ for
which
\begin{equation*}
P_{I}=P_{i_{0}}\cdots P_{i_{n-1}}=\mathbb{A}_{i_{0}}\cdots \mathbb{A}%
_{i_{n-1}},\mbox{ }\mbox{for all } I=\left( i_{0},\ldots ,i_{n-1}\right) \in
X_{n}
\end{equation*}%
where $\mathbb{A}_{k}(i,j)=\mathbb{A}\left( i,j\right) $ if $j=k$ and $%
\mathbb{A}_{k}(i,j)=0$ otherwise.

\item[(ii)] The push forward measure $\nu =\pi \mu $ is the linear
representable with respect to the triple $(x,Q,y)$, where $x=l,$ $y=\mathbf{1}%
_{n}$ and $Q$ is generated by $\left( Q_{j}\right) _{j\in \mathcal{A}%
(Y)}=\left( \mathbb{A}_{j}\right) _{j\in \mathcal{A}(Y)}$ for which $\mathbb{%
A}_{k}(u,v)=\mathbb{A}(u,v) $ if $\Pi (v)=k$ and $\mathbb{A}_{k}(u,v)=0$
otherwise.
\end{enumerate}
\end{proposition}

The following Proposition presents that the push forward measure of maximal
measure on $Z$ is the equilibrium measure with $N$. Recall that $%
A_{1},\ldots ,A_{\mathbf{n}}$ are induced from $A$ in Theorem \ref{Thm: 9},
define $\mathbb{A}_{j}= \left(stoch(A)\right)_{j}$ for $j= 1, \ldots,
\mathbf{n}$.

\begin{proposition}
\label{Prop: 1}Let $Z$ be a Markov Sierpi\'{n}ski carpet with $A$ being irreducible and the induced potential $N$ also being irreducible. Let $\mu _{A}$ be
the unique Markov measure of $A$, then $\nu =\bar{\nu}=\pi _{y}\mu _{A}.$
\end{proposition}

\begin{proof}
Since $\nu =\bar{\nu}$ from Corollary \ref{Cor: 1}, it suffices to show that $%
\nu =\pi _{y}\mu _{A}$. Since $\pi _{y}\mu _{A}$ is a linear representable
sofic measure by Proposition \ref{Prop: 3}. Let $\mathbb{A}= stoch\left(
A\right) $ and $l$ be the $1\times \mathbf{d}$ ($\textbf{d} = \textbf{m}\times \textbf{n}$) left eigenvector of $\mathbb{A%
}$ with respect to the maximal eigenvalue $1$. It follows from Proposition %
\ref{Prop: 3} that the triple $\left( l,\mathbb{A},\mathbf{1}_{\mathbf{d}%
}\right) $ defines a linear representable measure $\pi _{y}\mu _{A}$, where $%
\mathbb{A}=\left( \mathbb{A}_{j}\right) _{j=1}^{\mathbf{n}}$ with $\mathbb{A}%
_{j}(u,v)=\mathbb{A}\left( u,v\right) $ if $\pi _{y}\left( v\right) =j$, $%
\mathbb{A}_{j}\left( u,v\right) =0$ otherwise. That is,
\begin{equation*}
\pi _{y}\mu _{A}\left( J\right) =l\mathbb{A}_{j_{0}}\mathbb{A}_{j_{1}}\ldots
\mathbb{A}_{j_{n-1}}\mathbf{1}_{\mathbf{d}}\mbox{, for all }J=\left(
j_{0}\ldots j_{n-1}\right) \in Y_{n}\mbox{.}
\end{equation*}%
Under the same argument of the proof in Theorem \ref{Thm: 1} and $\log \rho
_{A}=P\left( \sigma _{Y},N\right) $. For $J=\left( j_{0}\ldots
j_{n-1}\right) \in Y_{n}$ we have%
\begin{eqnarray*}
\pi _{y}\mu _{A}\left( \left[ J\right] \right) &=&l\mathbb{A}_{j_{0}}\mathbb{%
A}_{j_{1}}\ldots \mathbb{A}_{j_{n-1}}\mathbf{1}_{\mathbf{d}} \\
&=&\rho _{A}^{-n}lD^{-1}A_{j_{0}}A_{j_{1}}\ldots A_{j_{n-1}}D\mathbf{1}_{%
\mathbf{d}} \\
&\leq &\max_{1\leq j\leq \mathbf{n}}\left\{ r^{-1}(j)\right\} \max_{1\leq
j\leq \mathbf{n}}\left\{ r(j)\right\} \rho _{A}^{-n}\mathbf{1}_{\mathbf{d}%
}^{t}A_{j_{0}}A_{j_{1}}\ldots A_{j_{n-1}}\mathbf{1}_{\mathbf{d}} \\
&=&c_{1}\rho _{A}^{-n}\left\Vert A_{j_{0}}A_{j_{1}}\ldots
A_{j_{n-1}}\right\Vert \\
&\leq &c_{2}\rho _{A}^{-n}\left\Vert N_{j_{0}j_{1}}N_{j_{1}j_{2}}\ldots
N_{j_{n-2}j_{n-1}}\right\Vert \\
&=&c_{2}\exp \left( -nP\left( \sigma_{Y},N\right) \right) \left\Vert
N_{j_{0}j_{1}}N_{j_{1}j_{2}}\ldots N_{j_{n-2}j_{n-1}}\right\Vert \\
&\leq &c_{3}\nu \left( \left[ J\right] \right) \mbox{ (Corollary \ref{Cor: 1}%
)}
\end{eqnarray*}%
Similarly, since $N$ is irreducible
\begin{eqnarray*}
\pi _{y}\mu _{A}\left( \left[ J\right] \right) &=&l\mathbb{A}_{j_{0}}\mathbb{%
A}_{j_{1}}\ldots \mathbb{A}_{j_{n-1}}\mathbf{1}_{d} \\
&=&\rho _{A}^{-n}lD^{-1}A_{j_{0}}A_{j_{1}}\ldots A_{j_{n-1}}D\mathbf{1}_{d}
\\
&\geq &\min_{1\leq j\leq \mathbf{n}}\left\{ r^{-1}(j)\right\} \min_{1\leq
j\leq \mathbf{n}}\left\{ r(j)\right\} \rho _{A}^{-n}\mathbf{1}%
_{d}^{t}A_{j_{0}}A_{j_{1}}\ldots A_{j_{n-1}}\mathbf{1}_{d} \\
&\geq &c_{5}\rho _{A}^{-n}\left\Vert N_{j_{0}j_{1}}N_{j_{1}j_{2}}\ldots
N_{j_{n-2}j_{n-1}}\right\Vert \\
&\geq &c_{6}\nu \left( \left[ J\right] \right)
\end{eqnarray*}%
Since $\pi _{y}\mu _{A}$ and $\nu $ are ergodic, $\pi _{y}\mu _{A}=\nu $%
. The proof is completed.
\end{proof}

\section{Proofs}

This section presents the detailed proofs for Theorem \ref{Thm: A_Rev: 1},
Theorem \ref{Thm: 8}, Theorem \ref{Thm: 3} and Theorem \ref{Thm: 6}.

\subsection{Proof of Theorem \protect\ref{Thm: A_Rev: 1}}

We first review some background knowledge of $\mathbf{a}$-weighted
thermodynamic formalism proposed by Barral and Feng \cite{BarralFeng2009} for $p$%
-specification shift space and by Feng \cite{Feng2011} for the weak $p$%
-specification case. For $\textbf{a} = (a, b)$, the
$\mathbf{a}$\emph{-weighted pressure} is defined as
follows:%
\begin{equation}
P^{\mathbf{a}}(\sigma _{X},\Phi )=\sup \left\{ \Phi _{\ast }(\eta
)+a h_{\eta }(\sigma _{X})+b h_{\pi \eta }(\sigma _{Y}):\eta \in
\mathcal{M}(X,\sigma _{X})\right\} .  \label{7}
\end{equation}%
Define the collection of equilibrium measures and $\mathbf{a}$\emph{%
-weighted equilibrium }as follows:
\begin{eqnarray*}
\mathcal{I}(\Phi ) &=\left\{ \mu \in \mathcal{M}(X,\sigma _{X}):\mu
\mbox{ is an equilibrium measure of }\Phi \right\} , \\
\mathcal{I}(\Phi ,\mathbf{a}) &=\left\{ \mu \in \mathcal{M}(X,\sigma
_{X}): \mu \mbox{ attain the supremum of (\ref{7})}\right\} .
\end{eqnarray*}%
Let $\pi :X\rightarrow Y$ be a factor, the \emph{conditional equilibrium
measure} $\mu \in \mathcal{M}(X,\sigma _{X})$ of $\Phi $ with respect to $%
\nu $ if $\pi \mu =\nu $, and $\mu $ satisfies the \emph{conditional
variational principle}, i.e.,
\begin{equation*}
\Phi _{\ast }(\mu )+h_{\mu }(\sigma _{X})=\sup \left\{ \Phi _{\ast }(\eta
)+h_{\eta }(\sigma _{X}):\eta \in \mathcal{M}(X,\sigma _{X}),\mbox{ }\pi
\eta =\nu \right\} .
\end{equation*}%
Denote by $\mathcal{I}_{\nu }\left( \Phi \right) $ the collection of all
conditional equilibrium measure of $\Phi $ with respect to $\nu $

\begin{theorem}[{\cite[Corollary 3.11]{Feng2011}},
{\cite[Theorem 1.1]{BarralFeng2009}}]
\label{Thm: 7}Let $\Phi =\left( \log \phi _{n}\right) _{n=1}^{\infty }$ be a
subadditive potential function on $X$. For all $J\in Y_{n}$, $n\in \mathbb{N}
$, define $\psi _{n}:X_{n}\rightarrow \mathbb{R}$ as follows%
\begin{equation*}
\psi _{n}(J)= \sum_{I\in X_{n}:\pi (I)=J}\phi \left( I\right).
\end{equation*}%
Let $\Psi =\left( \log \psi _{n}\right) _{n=1}^{\infty }$ be the collection
of $\psi _{n}$. Then

\begin{enumerate}
\item[(i)] $P^{\mathbf{a}}\left( \sigma _{X},\Phi \right)
=(a + b)P\left( \sigma _{Y},\left( \frac{a}{a+b}\right) \Psi \right) $
\item[(ii)] $\mu \in \mathcal{I}(\Phi ,\mathbf{a})$ if and only if $\mu
\circ \pi ^{-1}\in \mathcal{I}\left( \frac{a}{a+b}\Psi \right) $ and
$\mu \in \mathcal{I}_{\mu \circ \pi ^{-1}}\left( \frac{1}{a}\Phi \right)
$, where $\frac{a}{a+b}\Psi =\left( \left( \frac{a}{a+b}\right) \log
\left( \psi _{n}\right) _{n=1}^{\infty }\right) $ and $\frac{1}{a}\Phi
=\left( \left( \frac{1}{a}\right) \log \left( \phi _{n}\right) \right)
_{n=1}^{\infty }$.
\end{enumerate}
\end{theorem}

\begin{proof}[Proof of Theorem \protect\ref{Thm: A_Rev: 1}]
\textbf{Step 1.} ($(i):(a)\Rightarrow (b)$) We assume that $\mu $ is the
invariant measure of the full Hausdorff dimension, i.e., $\dim _{H}Z=\dim
_{H}\mu $. We claim that $\mu $ attains the supremum of (\ref{7}) with $%
\mathbf{a}=\left( \alpha ,1-\alpha \right) \in \mathbb{R}^{2}$ and zero
potential $\Phi $. Indeed, under the same argument of Theorem \ref{Thm: 1},
there exists a constant $c>0$ such that for all $y\in Y$, we have
\begin{equation}
c^{-1}\left\Vert N_{y|_{n}}\right\Vert \leq \psi _{n}(y|_{n})\leq
c\left\Vert N_{y|_{n}}\right\Vert \mbox{,}  \label{25}
\end{equation}%
where $y|_{n}=\left( y_{0},\ldots ,y_{n-1}\right) \in Y_{n}$. Combining
Theorem \ref{Thm: 1} and Theorem \ref{Thm: 7} we have
\begin{eqnarray}
\dim _{H}\mu &=&\dim _{H}Z=\frac{1}{\log \mathbf{n}}P(\sigma_{Y},\alpha N)\mbox{
(Theorem \ref{Thm: 1})}  \nonumber \\
&=&\frac{1}{\log \mathbf{n}}P(\sigma_{Y},\Psi )=\frac{1}{\log \mathbf{n}}P^{\mathbf{a}%
}(\sigma _{Z},\Phi )\mbox{.}  \label{41}
\end{eqnarray}%
Combining (\ref{41}) and the Ledrappier-Young formula for the Hausdorff
dimension of measure $\mu $ \cite{LedrappierYoung1985, LedrappierYoung1985a} yields
\begin{eqnarray*}
\frac{1}{\log \mathbf{n}}P^{\mathbf{a}}(\sigma _{Z},\Phi ) &=&\dim _{H}\mu \\
&=&\frac{h_{\mu }(\sigma _{Z})}{\log \mathbf{m}}+\left( \frac{1}{\log
\mathbf{n}}-\frac{1}{\log \mathbf{m}}\right) h_{\pi _{y}\mu }(\sigma _{Y}) \\
&=&\frac{h_{\mu }(\sigma _{Z})}{\log \mathbf{m}}+\frac{1}{\log \mathbf{n}}%
\left( 1-\frac{\log \mathbf{n}}{\log \mathbf{m}}\right) h_{\pi _{y}\mu
}(\sigma _{Y}) \\
&=&\frac{h_{\mu }(\sigma _{Z})}{\log \mathbf{m}}+\frac{1}{\log \mathbf{n}}%
\left( 1-\alpha \right) h_{\pi _{y}\mu }(\sigma _{Y})\mbox{.}
\end{eqnarray*}%
Therefore,
\begin{eqnarray*}
P^{\mathbf{a}}(\sigma _{Z},\Phi ) &=&\frac{\log \mathbf{n}}{\log \mathbf{m}}%
h_{\mu }(\sigma _{Z})+\left( 1-\alpha \right) h_{\pi _{y}\mu }(\sigma _{Y})
\\
&=&\alpha h_{\mu }(\sigma _{Z})+\left( 1-\alpha \right) h_{\pi _{y}\mu
}(\sigma _{Y})\mbox{.}
\end{eqnarray*}%
This shows that $\mu \in \mathcal{I}(\Phi ,\mathbf{a})$ with $\mathbf{a}%
=\left( \alpha ,1-\alpha \right) \in \mathbb{R}^{2}$.

\textbf{Step 2.} $\left( (i):(b)\Rightarrow (a)\right) $ It follows from the
Ledrappier-Young formula of measure $\mu $ and it is the $\mathbf{a}$%
-weighted equilibrium measure with $\mathbf{a}=\left( \alpha ,1-\alpha
\right) \in \mathbb{R}^{2}$ by Theorem \ref{Thm: 11}. Up to a minor
modification of Proposition 2.6 of \cite{BarralFeng2009}
\begin{equation}
N_{\ast }(\nu _{\alpha })+h_{\nu_{\alpha } }\left( \sigma _{Y}\right) =\sup
\left\{ h_{\eta }(\sigma _{Z}):\pi \eta =\nu _{\alpha }\right\}
\mbox{,} \label{42}
\end{equation}%
where
\[
N_{\ast }\left( \nu _{\alpha }\right) =\lim_{n\rightarrow \infty }\frac{1}{n}%
\int_{Y}\log \left\Vert N_{y|_{n}}\right\Vert d\nu _{\alpha
}(y)\mbox{.}
\]
It follows from the variational principle we have%
\begin{equation}
\left(\alpha N\right)_{\ast }(\nu _{\alpha })+h_{\nu _{\alpha }}\left( \sigma
_{Y}\right) =P(\sigma_{Y},\alpha N)\mbox{.}\label{44}
\end{equation}%
Combining Theorem \ref{Thm: 7} (ii), (\ref{42}), \eqref{44} and the Ledrappier-Young formula of measure $\mu $ obtains
\begin{eqnarray*}
&&\dim _{H}Z-\dim _{H}\mu \\
&=&\frac{1}{\log \mathbf{n}}\left( P(\sigma_{Y},\alpha N)-\alpha h_{\mu }(\sigma
_{Z})-(1-\alpha )h_{\nu _{\alpha }}(\sigma_{Y})\right) \\
&=&\frac{1}{\log \mathbf{n}}\left( P(\sigma_{Y},\alpha N)-h_{\nu _{\alpha }}(\sigma
_{Y})-\alpha N_{\ast }(\nu _{\alpha })\right) \\
&=&\frac{1}{\log \mathbf{n}}\left( P(\sigma_{Y},\alpha N)-\left( h_{\nu _{\alpha
}}(\sigma _{Y})+\left(\alpha N\right)_{\ast }(\nu _{\alpha })\right) \right) \\
&=&\frac{1}{\log \mathbf{n}}\left( P(\sigma_{Y},\alpha N)-P(\sigma_{Y},\alpha N)\right) =0%
\mbox{.}
\end{eqnarray*}%
This shows that $\mu $ is the invariant measure of the full Hausdorff dimension.
This completes the proof of (i).

\textbf{Step 3.} It remains to prove the dimension formula
(\ref{10}). Indeed, take $\mathbf{a}=(\alpha ,1-\alpha )\in
\mathbb{R}^{2}$ and $\nu _{\alpha }=\pi \mu \in \mathcal{I}(\alpha
N)$ and $\Phi $ is zero potential.
It follows from Theorem \ref{Thm: A_Rev: 1}-(i) and the definition of $%
\mathbf{a}$-weighted pressure (\ref{7}) we obtain that
\begin{eqnarray}
P^{\mathbf{a}}(\sigma _{Z},\Phi ) &=&\sup \left\{ \alpha h_{\eta }(\sigma _{%
Z})+(1-\alpha )h_{\pi _{y}\eta }(\sigma _{\mathbf{Y}}):\eta \in
\mathcal{M}(Z,\sigma _{Z})\right\} .  \nonumber \\
&=&\alpha h_{\mu }(\sigma _{Z})+(1-\alpha )h_{\pi _{y}\mu }(\sigma _{Y})
\label{26}
\end{eqnarray}%
Combining (\ref{26}) and the fact of $h_{top}(\sigma _{Z})=\log \rho _{A}$
we have
\begin{eqnarray*}
\dim _{H}Z &=&\frac{1}{\log \mathbf{n}}P(\sigma_{Y},\alpha N)=\frac{1}{\log \mathbf{n}%
}(h_{\nu _{\alpha }}(\sigma _{Y})+\left(\alpha N\right)_{\ast }(\nu _{\alpha })) \\
&=&\frac{1}{\log \mathbf{n}}(h_{\bar{\nu}_{\alpha }}(\sigma
_{Y})+\alpha
N_{\ast }(\bar{\nu}_{\alpha }))\mbox{ (Proposition \ref{Prop: 1})} \\
&=&\frac{1}{\log \mathbf{n}}\left( \alpha\log \rho
_{A}+h_{\bar{\nu}_{\alpha
}}(\sigma _{Y})+\left(\alpha \bar{N}\right)_{\ast }(\bar{\nu}_{\alpha })\right) \mbox{ (}%
\bar{N}\mbox{ is normalized)} \\
&=&\frac{1}{\log \mathbf{n}}\left( \alpha h_{top}(\sigma _{Z})+P(\sigma
_{Y},\alpha \bar{N})\right) \\
&=&\frac{\alpha h_{top}(Z)}{\log \mathbf{n}}+\frac{P(\sigma _{Y},\alpha \bar{N})}{%
\log \mathbf{n}} \\
&=&\frac{h_{top}(Z)}{\log \mathbf{m}}+\frac{P(\sigma _{Y},\alpha
\bar{N})}{\log \mathbf{n}}\mbox{.}
\end{eqnarray*}%
This establishes the formula (\ref{10}).
\end{proof}

\subsection{Proof of Theorem \protect\ref{Thm: 8}}

For the proof of Theorem \ref{Thm: 8} we give some useful lemmas first.

\begin{lemma}
\label{Lma: 7}Let $A\in \mathbb{R}^{n\times n}$, $A$ be irreducible
and $L=\left( \left[ L\left( i\right) \right] _{i=1}^{n}\right)
^{t}$, $R=\left[ R\left( i\right) \right] _{i=1}^{n}$ be the left
and right eigenvector of $A$ corresponding to the
maximal eigenvalue $\rho _{A}$. If rank$\left( A\right) =1$, then $L=\left( %
[ C\left( j\right) ] _{j=1}^{n}\right) ^{t}$ and
$R=\left[ D\left( i\right) \right] _{i=1}^{n}$, where
\[
C\left( j\right) =\left\{
                    \begin{array}{ll}
                      A\left( i,j\right) /A\left( i,1\right) , & \hbox{ if $A\left( i,1\right)\neq 0$;} \\
                      0, & \hbox{otherwise.}
                    \end{array}
                  \right.
\]
and
\[
D\left( i\right) =\left\{
                    \begin{array}{ll}
                      A\left( i,j\right) /A\left( 1, j\right) , & \hbox{ if $A\left( 1, j\right)\neq 0$;} \\
                      0, & \hbox{otherwise.}
                    \end{array}
                  \right.
\]
\end{lemma}

\begin{proof}
Without loss of generality, we may assume that $L(1)=1$. Since
rank$(A)=1$, $C(j)$ is well-defined for all $j =1, \ldots,
\textbf{n}$. It follows from the fact that $L$ is the left
eigenvector of $A$ with respect to the eigenvalue $\rho _{A}$, we
have for all $j=2,\ldots ,n.$
\begin{equation*}
\sum_{i=1}^{n}L\left( i\right) A\left( i,j\right) =\sum_{i=1}^{n}L\left(
i\right) C\left( j\right) A\left( i,1\right) =C\left( j\right) \rho
_{A}L\left( 1\right) \mbox{.}
\end{equation*}%
Hence $L\left( j\right) =C\left( j\right) L\left( 1\right) $ for all $%
j=1,\ldots ,n$ and $L=\left([ C\left( j\right) ]
_{j=1}^{n}\right) ^{t}$. It follows from the same argument that we have $R=\left[
D\left( i\right) \right] _{i=1}^{n}$. This completes the proof.
\end{proof}

\begin{lemma}
\label{Lma: 6}Let $\widehat{L}=\left( \left[ \widehat{L}\left( k\right) %
\right] _{k=1}^{\mathbf{n}}\right) ^{t}$, where $\widehat{L}\left( k\right) =\left( %
\left[ L\left( d_{1},k\right) \right] _{d_{1}=1}^{\mathbf{m}}\right) ^{t}$
and $B_{k}=PA_{k}P^{-1}$ for all $k=1,\ldots ,\mathbf{n}$, where $%
B_{k}^{\left( i,j\right) }=N_{ij}$ if $j=k,$ and $B_{k}^{\left( i,j\right)
}=0_{\mathbf{m}\times \mathbf{m}}$, otherwise. Then for all $k=1,\ldots ,%
\mathbf{n}$,
\begin{equation*}
\widehat{L}B_{k}=\left[ 0_{\mathbf{m}},0_{\mathbf{m}},\ldots ,\rho _{A}%
\widehat{L}\left( k\right) ,0_{\mathbf{m}},\ldots ,0_{\mathbf{m}}\right]
\mbox{,}
\end{equation*}
where $0_{\mathbf{m}}\in \mathbb{R}^{1\times \mathbf{m}}$ with all entries
being $0^{\prime }s$.
\end{lemma}

\begin{proof}
Since $LA=\rho _{A}L$, we have $\sum_{d\in D}L\left( d\right) A\left(
d,d^{\prime }\right) =\rho _{A}L\left( d^{\prime }\right) $ for all $%
d^{\prime }\in D$. For each $1\leq k\leq \mathbf{n}$, construct $A_{k}$ as
follows:%
\begin{equation*}
A_{k}\left( d,d^{\prime }\right) =\left\{
\begin{array}{ccc}
A\left( d,d^{\prime }\right), &  & \mbox{if }d^{\prime }=\left( d_{1}^{\prime
},d_{2}^{\prime }\right) \mbox{ with }d_{2}^{\prime }=k; \\
0, &  & \mbox{otherwise.}%
\end{array}%
\right.
\end{equation*}%
We obtain
\begin{equation*}
\sum_{d\in D}L\left( d\right) A_{k}\left( d,d^{\prime }\right) =\left\{
\begin{array}{ccc}
\rho _{A}L\left( d^{\prime }\right), &  & \mbox{if }d^{\prime }=\left(
d_{1}^{\prime },d_{2}^{\prime }\right) \mbox{ with }d_{2}^{\prime }=k; \\
0, &  & \mbox{otherwise.}%
\end{array}%
\right.
\end{equation*}%
On the other hand, since $\widehat{L} = L P^{-1}$ and $LA_{k}=L\left( P^{-1}B_{k}P\right) $ for all $%
k=1,\ldots ,\mathbf{n}$,
\begin{equation*}
\widehat{L}B_k = LA_{k}P^{-1}=\rho _{A}\left[ 0_{\mathbf{m}},\ldots ,0_{\mathbf{m}},\widehat{L%
}\left( k\right) ,0_{\mathbf{m}},\ldots ,0_{\mathbf{m}}\right] \mbox{. }
\end{equation*}%
The proof is thus completed.
\end{proof}

Notably, if we let $\widehat{R}=PR=\left[ \widehat{R}\left( k\right) \right]
_{k=1}^{\mathbf{n}}$, where $\widehat{R}\left( k\right) =\left[ R\left(
d_{1},k\right) \right] _{d_{1}=1}^{\mathbf{m}}$, we also have
\begin{equation*}
B_{k}\widehat{R}=\left( \left[ 0_{\mathbf{m}},\ldots ,0_{\mathbf{m}},\rho
_{A}\widehat{R}\left( k\right) ,0_{\mathbf{m}},\ldots ,0_{\mathbf{m}}\right]
\right) ^{t}\mbox{ for all }k=1,\ldots ,\mathbf{n.}
\end{equation*}

\begin{proof}[Proof of Theorem \protect\ref{Thm: 8}]
For clarity, we prove theorem \ref{Thm: 8} for the cases of $k=1$ and $k\geq
2$.

\textbf{Step 1 }$\mathbf{(}k=1)$.

\textbf{1}.$(\Rightarrow)$ We claim that if $Z$ is a Markov
Sierpi\'{n}ski carpet, then the induced potential $N$ satisfies the
Markov condition of order $1$, that is, there exists $V_{1},\ldots ,V_{%
\mathbf{n}}$ such that $V_{i}N_{ij}\subseteq V_{j}$ for all $1\leq
i,j\leq \mathbf{n}$.

Since $Z$ is a Markov Sierpi\'{n}ski carpet, we obtain rank$%
\left( A\right) =1$. Hence, for all $1\leq j\leq \mathbf{m}, 1\leq
k\leq \mathbf{n}$ and $d\in D$,
\begin{equation}
C\left( j, k\right)=\left\{
                      \begin{array}{ll}
                        A\left( d,\left( j, k\right) \right)/ A\left( d,\left( 1,k\right)\right), & \hbox{ if $A\left( d,\left( 1,k\right)\right)\neq 0$;} \\
                        0, & \hbox{otherwise}
                      \end{array}
                    \right.  \label{21}
\end{equation}%
is well-defined.
Let $L=\left( \left[ L\left( d\right) \right] _{d\in D}\right) ^{t}$ be the
left eigenvector of $A$ corresponding to the maximal eigenvalue $\rho _{A}$, we
have $\widehat{L}:=LP^{-1}=\left( \left[ \widehat{L}\left( k\right) \right]
_{k=1}^{\mathbf{n}}\right) ^{t}$, where $\widehat{L}\left( k\right) =\left( %
\left[ L\left( d_{1},k\right) \right] _{d_{1}=1}^{\mathbf{m}}\right) ^{t}$. It follows from (\ref{21}) and Lemma \ref{Lma: 7}, we obtain
\begin{equation}
\widehat{L}\left( k\right) =\left( \left[ C\left( j, k\right) \right] _{j=1}^{%
\mathbf{m}}\right) ^{t}.\label{1}
\end{equation}

Choose $V_{k}=\widehat{L}\left( k\right) $ for $1\leq k\leq \mathbf{n}$.  Combining the fact that $N_{ij}=\left[ A\left( \left( p,i\right) ,\left(
q,j\right) \right) \right] _{p,q=1}^{\mathbf{m}}$ and  (\ref{21}), yields
\begin{equation}
A\left( \left( p,i\right)
,\left( q,j\right) \right) =C\left( q, j\right) A\left( \left(
p,i\right) ,\left( 1,j\right) \right) \label{2}
\end{equation}
for all $1\leq q\leq \mathbf{m}$, $1\leq j\leq \mathbf{n%
}$  and $\left( p,i\right) \in D$.

According to (\ref{1}) and (\ref{2}), for any $i,j=1,\ldots ,%
\mathbf{n}$, we get
\begin{eqnarray*}
V_{i}N_{ij} &=&\left[ C(1, i),C\left( 2,i\right) ,\ldots ,C\left(\mathbf{m}, i%
\right) \right] \times \\
&&\left[
\begin{array}{ccc}
A\left( \left( 1,i\right) ,\left( 1,j\right) \right) & \cdots & A\left(
\left( 1,i\right) ,\left( \mathbf{m},j\right) \right) \\
\vdots & \ddots & \vdots \\
A\left( \left( \mathbf{m},i\right) ,\left( 1,j\right) \right) & \cdots &
A\left( \left( \mathbf{m},i\right) ,\left( \mathbf{m},j\right) \right)%
\end{array}%
\right] \\
&=&\sum_{p=1}^{\mathbf{m}}C\left( p, i\right) A\left( \left(
p,i\right)
,\left( 1,j\right) \right) \left( \left[ C\left( q, j\right) \right] _{q=1}^{%
\mathbf{m}}\right) ^{t}=m\left( i,j\right) V_{j}\mbox{,}
\end{eqnarray*}%
where $m\left( i,j\right) =\sum_{p=1}^{\mathbf{m}}C\left( p,
i\right) A\left( \left( p,i\right) ,\left( 1,j\right) \right) $.
Therefore, $N$ satisfies the Markov condition of order $1$.

\textbf{2}.$(\Leftarrow)$ For the converse, we show that if $N$ satisfies the Markov condition of order $1$, then $%
\nu $ is a $1$-step Markov measure, i.e.,
\begin{equation*}
\frac{\nu \left( \left[ j_{0}j_{1}\right] \right) }{\nu \left( \left[ j_0%
\right] \right) }=\frac{\nu \left( \left[ j_{-k}\ldots j_{0}j_{1}\right]
\right) }{\nu \left( \left[ j_{-k}\ldots j_{0}\right] \right) },\mbox{ for all } k> 0.
\end{equation*}%
Assume $N$ satisfies the Markov condition of order $1$, it follows from the same argument
as above, that we have
\begin{equation}
\widehat{L}\left( i\right) N_{ij}=m\left( i,j\right)
\widehat{L}\left( j\right) , \hbox{ for all } 1\leq i,j\leq \mathbf{n}.\label{31}
\end{equation}
Since for each $J=\left[ j_{0}\ldots j_{n-1}\right] \in Y_{n}$, $\nu \left( J\right)
=\left( 1/\rho_A\right) ^{n}LA_{j_{0}}A_{j_{1}}\cdots A_{j_{n-1}}R$. Hence,
it follows from Lemma \ref{Lma: 6} and (\ref{31})
\begin{eqnarray*}
&&\frac{\nu \left( \left[ j_{0}j_{1}\right] \right) }{\nu \left( \left[ j_{0}%
\right] \right) }=\left(\frac{1}{\rho _{A}\left[ 0_{\mathbf{m}},\ldots ,0_{%
\mathbf{m}},\rho_A\widehat{L}\left( j_{0}\right) ,0_{\mathbf{m}},\ldots ,0_{%
\mathbf{m}}\right] \widehat{R}}\right)\times \\
&&\stackrel{\mbox{the only non-zero part is in }j_{0}\mbox{-th coordinate}}{%
\overbrace{\left[ 0_{\mathbf{m}},\ldots ,0_{\mathbf{m}},\rho_A\widehat{L}%
\left( j_{0}\right) ,0_{\mathbf{m}},\ldots ,0_{\mathbf{m}}\right] }}\times\\
&&\stackrel{\mbox{the only non-zero column is in }j_{1}\mbox{-th coordinate}}{%
\overbrace{\left[
\begin{array}{ccccc}
0_{\mathbf{m}} & \cdots & \ N_{1j_{1}} & \cdots & 0_{\mathbf{m}} \\
\vdots & \vdots & \vdots & \vdots & \vdots \\
0_{\mathbf{m}} & \cdots & N_{nj_{1}} & \cdots & 0_{\mathbf{m}}%
\end{array}%
\right] }}\times \widehat{R} \\
&=&\frac{m\left( j_{0},j_{1}\right) \widehat{L}\left(j_1\right)\widehat{R}%
\left( j_{1}\right) }{\rho_{A}\widehat{L}\left( j_{0}\right) \widehat{R}%
\left( j_{0}\right) }\mbox{.}
\end{eqnarray*}%
On the other hand, for any $k > 0$,%
\begin{eqnarray}
&&\frac{\nu \left( \left[ j_{-k}j_{-k+1}\cdots j_{0}j_{1}\right] \right) }{%
\nu \left( \left[ j_{-k}\cdots j_{0}\right] \right) }=\left(\frac{1}{\rho_A%
\widehat{L}B_{j_{-k}}\cdots B_{j_{0}}\widehat{R}}\right)\times \nonumber\\
&&\stackrel{\mbox{the non-zero element is in the }j_{-k+1}\mbox{-th coordinate%
}}{\overbrace{\left[ 0_{\mathbf{m}},\ldots ,0_{\mathbf{m}},\widehat{L}\left(
j_{-k}\right) N_{j_{-k}j_{-k+1}},0_{\mathbf{m}},\ldots ,0_{\mathbf{m}}\right]
}}B_{j_{-k+2}}\cdots B_{j_{0}}B_{j_{1}}\widehat{R} \nonumber\\
&=&\left(\frac{1}{\rho_A\widehat{L}B_{j_{-k}}\cdots B_{j_{0}}\widehat{R}}%
\right)m\left( j_{-k},j_{-k+1}\right) \times \nonumber\\
&&\left[ 0_{\mathbf{m}},\ldots ,0_{\mathbf{m}},\widehat{L}\left(
j_{-k+1}\right) ,0_{\mathbf{m}},\ldots ,0_{\mathbf{m}}\right]
B_{j_{-k+2}}\cdots B_{j_{0}}B_{j_{1}}\widehat{R}.\label{3}
\end{eqnarray}%
Continuing the process as (\ref{3}), we have
\begin{equation*}
\frac{\nu \left( \left[ j_{-k}j_{-k+1}\cdots j_{0}j_{1}\right] \right) }{\nu
\left( \left[ j_{-k}\cdots j_{0}\right] \right) }=\frac{m\left(
j_0,j_1\right) \widehat{L}\left( j_1\right) \widehat{R}\left( j_1\right) }{%
\rho_A\widehat{L}\left( j_{0}\right) \widehat{R}\left( j_{0}\right) }=\frac{%
\nu \left( \left[ j_{0}j_{1}\right] \right) }{\nu \left( \left[ j_{0}\right]
\right) }
\end{equation*}for all $k > 0.$
Hence, $\nu $ is a $1$-step Markov measure and the proof is thus completed.

\textbf{Step 2} ($k\geq 2$). For this proof, we recall some setting first, since $Z=Z_{\left(
\mathbf{m},\mathbf{n}\right) }\left( A\right) $ is a Markov Sierpi\'{n}ski
carpet, denote $Z^{\left[ k\right] }=Z_{\left(\mathbf{m}^k,\mathbf{n}^k%
\right) }\left( A^{\left[ k\right] }\right) $, where $A^{\left[
k\right] }\in \mathbb{R}^{\textbf{d}^k \times \textbf{d}^k}$
($\textbf{d} = \textbf{m}\times \textbf{n}$) is the $k$-th higher
block transition matrix from $A$ which is indexed by $D^{[k]}$ and
define as $P^{[k]} = P^{[k]}_{\left( \mathbf{m},\mathbf{n}\right)}$.
Let $L^{[k]} = \left(\left[L^{[k]}(d^{[k]})\right]_{d^{[k]}\in
D^{[k]}}\right)^t$ be the left eigenvector of $A^{[k]}$
corresponding to the maximal eigenvalue $\rho_{A^{[k]}} = \rho_A$.
Set $B^{\left[ k\right] }=P^{[k]}A^{\left[ k\right]
}\left(P^{[k]}\right)^{-1}$. For $J\in Y_k$, we construct
$A^{[k]}_J$ as follows:
\begin{equation}
\label{56}
A^{[k]}_J\left(d^{[k]}, d'^{[k]}\right) = \left\{
                                            \begin{array}{ll}
                                              A^{[k]}\left(d^{[k]}, d'^{[k]}\right), & \hbox{ if  $d^{[k]} = (d^{[k]}_1, d^{[k]}_2)$ with $d^{[k]}_2 = J$}; \\
                                              0, & \hbox{ otherwise}
                                            \end{array}
                                          \right.
\end{equation}
and define an ordered set
\[
\Gamma_J := \left\{\left(d^{[k]}\right)^{(r)} :
\left(d^{[k]}\right)^{(r)}=\left(\left(d^{[k]}\right)^{(r)}_1,
\left(d^{[k]}\right)^{(r)}_2\right) \hbox{ with }
\left(d^{[k]}\right)^{(r)}_2 = J\right\}.
\]
It can be easily checked that $\#\Gamma_J = \mathbf{m}^k$ for any $J\in Y_k.$

We present two lemmas which are analogous to Lemma \ref{Lma: 1} and  \ref{Lma: 6} for $k=1$.

\begin{lemma}[Lemma \ref{Lma: 1} for $k \geq 2$]
Let $Z=Z_{\left( \mathbf{m},\mathbf{n}\right)}(A)$ be a Markov
Sierpi\'{n}ski carpet with $A\in \mathbb{R}%
^{\mathbf{d}\times \mathbf{d}}$ $(\mathbf{d}=\mathbf{m}\times \mathbf{n})$ be irreducible. For $k\geq 2$%
, $A^{\left[ k\right] }$ is the $k$-th higher block transition matrix which
is indexed by $D^{[k]}$. For all $J \in Y_{k}$, let $A_{J}^{\left[ k\right] }$
be  as defined in (\ref{56}) and $B_{J}^{[k]}=P^{[k]}A^{[k]}_{J}\left(P^{[k]}\right)^{-1}$. Then
\begin{equation*}
\left( B_{J}^{\left[ k\right] }\right) ^{\left( J^{\prime }\left(
0,k-1\right) , J^{\prime }\left( 1,k\right) \right) } = \left\{
\begin{array}{ccc}
N_{J^{\prime }\left( 0,k\right) }, &  & \mbox{if }J^{\prime }\left(
1,k\right) =J\mbox{;} \\
0_{\mathbf{m}^{k}\times \mathbf{m}^{k}}, &  & \mbox{otherwise.}%
\end{array}%
\right.
\end{equation*}
\end{lemma}

\begin{lemma}[Lemma \ref{Lma: 6} for $k \geq 2$]
\label{Lma: 8}Let $\widehat{L}^{\left[ k\right] }:=L^{\left[ k\right]
}\left(P^{[k]}\right)^{-1}=\left( \left[ \widehat{L}^{\left[ k\right] }\left( J_{i}\right) %
\right] _{i=1}^{a_{k}}\right)^t $, where
\begin{equation*}
\widehat{L}^{\left[ k\right] }\left( J_{i}\right) =\left( \left[ L^{\left[ k%
\right] }\left( \left(d^{[k]}\right)^{(r)}\right) \right] _{r=1}^{\mathbf{m}%
^{k}}\right) ^{t}\mbox{,}
\end{equation*}
for all $\left(d^{[k]}\right)^{(r)}\in \Gamma_{J_{i}}$, $r = 1,\ldots, \mathbf{m}^k$ and $%
B_{J_{i}}^{\left[ k\right] }=P_kA_{J_{i}}^{\left[ k\right]
}P_k^{-1}$ for all $i=1,\ldots ,a_{k}$. Then
\begin{equation*}
\widehat{L}^{\left[ k\right] }B_{J_i}^{\left[ k\right] }=\left[ 0_{\mathbf{m}%
^{k}},\ldots ,0_{\mathbf{m}^{k}},\rho _{A}\widehat{L}\left( J_{i}\right) ,0_{%
\mathbf{m}^{k}},\ldots ,0_{\mathbf{m}^{k}}\right] \mbox{ for all }
i=1,\ldots ,a_{k}.
\end{equation*}
\end{lemma}
We continue the proof of \textbf{Step 2} and divide it into two small parts.
\textbf{1}.$(\Rightarrow)$ We claim that if $Z$ is a Markov
Sierpi\'{n}ski carpet, then $N$ satisfies the Markov condition of order k, that is, there exists $\left\{ V_{J}\right\} _{J\in Y_{k}}$ such that $%
V_{J\left( 0,k-1\right) }N^{[k]}_{J\left( 0,k\right) }\subseteq
V_{J\left( 1,k\right) }$ for all $J\left( 0,k-1\right) ,$ $J\left(
1,k\right) \in Y_{k}$.
Indeed, since $Z^{\left[ k\right] }$ is also
a Markov Sierpi\'{n}ski carpet, we have $rank(A^{[k]}) = 1$. Hence,
for all $J\in Y_k$,
\begin{equation}
\label{22} A^{[k]}\left(d^{[k]}, (d'^{[k]})^{(r)}\right) = C(J,
r)A^{[k]}\left(d^{[k]}, (d'^{[k]})^{(1)}\right),
\end{equation}
where $(d'^{[k]})^{(r)}\in \Gamma_J, r = 1, \ldots, \textbf{m}^k$
and $d^{[k]}\in D^{[k]}.$ let $L^{[k]} =
\left(\left[L^{[k]}(d^{[k]})\right]_{d^{[k]}\in
D^{[k]}}\right)^t$ be the left eigenvector of $A^{\left[ k%
\right] }$ corresponding to the maximal eigenvalue $\rho _{A^{\left[
k\right] }}=\rho _{A}$. This
implies $\widehat{L}^{\left[ k\right] }:=L^{\left[ k\right] }P^{-1}=\left( %
\left[ \widehat{L}^{\left[ k\right] }\left( J\right) \right] _{J\in
Y_{k}}\right) ^{t}$, where $\widehat{L}^{\left[ k\right] }\left( J\right)
=\left( \left[ L^{\left[ k%
\right] }\left( \left(d^{[k]}\right)^{(r)}\right) \right] _{r=1}^{\mathbf{m}%
^{k}}\right) ^{t}$, for all $J\in Y_{k}$,
$\left(d^{[k]}\right)^{(r)}\in \Gamma_{J}$. It follows from Lemma
\ref{Lma: 7} and (\ref{22}) we have
\begin{equation}
\widehat{L}^{\left[ k\right] }\left( J\right) =\left( \left[ C\left(
J,r\right) \right] _{r=1}^{\mathbf{m}^{k}}\right) ^{t}\hbox{ for all
} J\in Y_{k}.\label{4}
\end{equation}

Taking $V_{J}:= \widehat{L}^{\left[ k\right] }\left( J\right) $ for
all $J\in Y_{k}$. For any $J=J\left( 0,k\right) \in Y_{k+1}$ we have
\[
N_{J} = N_{J\left( 0,k\right)} = \left[ A^{[k] }\left( (d^{[k]})^{(p)}, (d'^{[k]})^{(q)}\right) \right] _{p,q=1}^{\mathbf{%
m}^{k}},
\]
where $(d^{[k]})^{(p)}\in \Gamma_{J\left( 0,k-1\right) }$, $%
(d'^{[k]})^{(q)}\in \Gamma_{J\left( 1,k\right) }$. It also follows
from (\ref{22}) we obtain that
\begin{equation}
A^{\left[ k%
\right] }\left((d^{[k]})^{(p)}, (d'^{[k]})^{(q)}\right) = C\left(
J\left( 1, k\right) ,q\right) A^{\left[ k\right]
}\left((d^{[k]})^{(p)}, (d'^{[k]})^{(1)}\right),\label{8}
\end{equation}
for all $(d^{[k]})^{(p)}\in \Gamma_{J\left( 0,k-1\right) }$, $%
(d'^{[k]})^{(q)}\in \Gamma_{J\left( 1,k\right) }$ and $p, q=1,\ldots
,\mathbf{m}^{k}$. We conclude from (\ref{4}) and (\ref{8}) that for
any $J
= J\left( 0,k\right) \in Y_{k+1}$ we have%
\begin{eqnarray*}
&&V_{J\left( 0,k-1\right) }N_{J\left( 0,k\right) }=\left[ 1,C\left( J\left(
0,k-1\right) ,2\right) ,\ldots ,C\left( J\left( 0,k-1\right) ,\mathbf{m}%
^{k}\right) \right] \times \\
&&\left[
\begin{array}{ccc}
A^{[k] }\left( (d^{[k]})^{(1)}, (d'^{[k]})^{(1)}\right) & \cdots & A^{[k] }\left( (d^{[k]})^{(1)}, (d'^{[k]})^{(\mathbf{m}^k)}\right) \\
\vdots & \ddots & \vdots \\
A^{[k] }\left( (d^{[k]})^{(1)}, (d'^{[k]})^{(\mathbf{m}^k)}\right) &
\cdots &
A^{[k] }\left( (d^{[k]})^{(\mathbf{m}^k)}, (d'^{[k]})^{(\mathbf{m}^k)}\right)%
\end{array}%
\right] \\
&=&m\left( J\left( 0,k-1\right) ,J\left( 1,k\right) \right) V_{J\left(
1,k\right) }\mbox{,}
\end{eqnarray*}%
where $m\left( J\left( 0,k-1\right) ,J\left( 1,k\right) \right)
:=\sum_{p=1}^{\mathbf{m}^{k}}C\left( J\left( 0,k-1\right) ,p\right)
A^{\left[ k\right] }\left( (d^{[k]})^{(p)}, (d'^{[k]})^{(1)}\right)
$.

Hence, there exists $\left\{ V_{J}\right\} _{J\in Y_{k}}$ such that $%
V_{J\left( 0,k-1\right) }N^{[k]}_{J\left( 0,k\right) }\subseteq
V_{J\left(
1,k\right) }$ for all $J\left( 0,k-1\right) ,$ $J\left( 1,k\right) \in Y_{k}$%
, i.e., $N$ satisfies the Markov condition from left of order $k$.

\textbf{2}.$(\Leftarrow)$ we claim that if $N$ satisfies the Markov
condition from left of order $k$, then $\nu $ is a $k$-step Markov
measure on $Y$, i.e.,
\begin{equation*}
\frac{\nu \left( J\left( 0,k\right) \right) }{\nu \left( J\left(
0,k-1\right) \right) }=\frac{\nu \left( J\left( -n,k\right) \right) }{\nu
\left( J\left( -n,k-1\right) \right) }\mbox{, for all } J\left( -n,k\right)
\in Y_{n+k+1},\mbox{ }n \in\mathbb{N}.
\end{equation*}%
Assume $N$ satisfies the Markov condition from left of order $k$, it
follows from the same argument as above, we have $V_{J}:=
\widehat{L}^{\left[ k\right] }\left( J\right) $ for all $J\in Y_{k}$
and for all $J\left( 0,k\right) \in Y_{k+1}$,
\begin{equation}
\widehat{L}^{\left[ k\right] }\left( J\left( 0,k-1\right) \right) N_{J\left(
0,k\right) }=m\left( J\left( 0,k-1\right) ,J\left( 1,k\right) \right)
\widehat{L}^{\left[ k\right] }\left( J\left( 1,k\right) \right) .\label{29}
\end{equation}
On the other hand, $Z = Z_{\left( \mathbf{m},\mathbf{n}\right)
}\left( A\right) =Z_{\left( \mathbf{m}^k,\mathbf{n}^k\right) }\left(
A^{\left[ k\right] }\right) $ for $k\geq 2$, hence for any
$J=J\left( 0,m-1\right) \in Y_{m}$
with $m\geq k$%
\begin{eqnarray*}
\nu \left( J\left( 0,m-1\right) \right) &=&\left( 1/\rho _{A}\right)
^{m-k+1}L^{\left[ k\right] }A_{J\left( 0,k-1\right) }^{\left[ k\right]
}A_{J\left( 1,k\right) }^{\left[ k\right] }\cdots A_{J\left( m-k,m-1\right)
}^{\left[ k\right] }R^{\left[ k\right] } \\
&=&\left( 1/\rho _{A}\right) ^{m-k+1}L^{\left[ k\right] }\left(
\left(P^{[k]}\right)^{-1}B_{J\left( 0,k+1\right) }^{\left[ k\right]
}P^{[k]}\right)\times \\
&& \cdots \times \left(
\left(P^{[k]}\right)^{-1}B_{J\left( m-k,m-1\right) }^{\left[ k\right] }P^{[k]}\right) \\
&=&\left( 1/\rho _{A}\right) ^{m-k+1}\widehat{L}^{\left[ k\right]
}B_{J\left( 0,k+1\right) }^{\left[ k\right] }\cdots B_{J\left(
m-k,m-1\right) }^{\left[ k\right] }\widehat{R}^{\left[ k\right]
}\mbox{.}
\end{eqnarray*}%
Combining the above computation, (\ref{29}) and Lemma \ref{Lma: 8} yields
\begin{eqnarray}
\nu \left( J\left( 0,m-1\right) \right) &=&\left( 1/\rho _{A}\right)
^{m-k}m\left( J\left( 0,k-1\right) ,J\left( 1,k\right) \right) \times \nonumber\\
&&\left[ 0_{\mathbf{m}^{k}},\ldots ,0_{\mathbf{m}^{k}},\widehat{L}^{\left[ k%
\right] }\left( J\left( 1,k\right) \right) ,0_{\mathbf{m}^{k}},\ldots ,0_{%
\mathbf{m}^{k}}\right] \times \nonumber\\
&&B_{J\left( 2,k+1\right) }^{\left[ k\right] }\cdots
B_{J\left( m-k,m-1\right) }^{\left[ k\right] }\widehat{R}^{\left[ k\right] }\label{24}
\end{eqnarray}%
Continuing the same process as (\ref{24}), we obtain
\begin{eqnarray*}
&&\nu \left( J\left( 0,m-1\right) \right) \\
&=&\left( 1/\rho _{A}\right) ^{m-k}m\left( J\left( 0,k-1\right) ,J\left(
1,k\right) \right) \cdots m\left( J\left( m-k-1,m-2\right) ,J\left(
m-k,m-1\right) \right) \\
&&\widehat{L}^{\left[ k\right] }\left( J\left( m-k,m-1\right) \right)
\widehat{R}\left( J\left( m-k,m-1\right) \right)
\end{eqnarray*}%
Hence
\begin{equation*}
\frac{\nu \left( J\left( 0,k\right) \right) }{\nu \left( J\left(
0,k-1\right) \right) }=\frac{m\left( J\left( 0,k-1\right) ,J\left(
1,k\right) \right) \widehat{L}^{\left[ k\right] }\left( J\left( 1,k\right)
\right) \widehat{R}^{\left[ k\right] }\left( J\left( 1,k\right) \right) }{%
\rho _{A}\widehat{L}^{\left[ k\right] }\left( J\left( 0,k-1\right) \right)
\widehat{R}^{\left[ k\right] }\left( J\left( 0,k-1\right) \right) }.
\end{equation*}%
On the other hand, for any $n > 0$,%
\begin{eqnarray*}
&&\frac{\nu \left( J\left( -n,k\right) \right) }{\nu \left( J\left(
-n,k-1\right) \right) } \\
&=&\frac{\rho _{A}^{-\left(n+1\right) }\prod_{i=0}^{n}m\left( J\left(
-n+i,-n+i+k-1\right) ,J\left( -n+i+1,-n+i+k\right) \right) }{\rho
_{A}^{-n}\prod_{i=0}^{n-1}m\left( J\left( -n+i,-n+i+k-1\right) ,J\left(
-n+i+1,-n+i+k\right) \right) }\times \\
&&\frac{\widehat{L}^{\left[ k\right] }\left( J\left( 1,k\right) \right)
\widehat{R}^{\left[ k\right] }\left( J\left( 1,k\right) \right) }{\widehat{L}%
^{\left[ k\right] }\left( J\left( 0,k-1\right) \right) \widehat{R}^{\left[ k%
\right] }\left( J\left( 0,k-1\right) \right) } \\
&=&\frac{m\left( J\left( 0,k-1\right) ,J\left( 1,k\right) \right) \widehat{L}%
^{\left[ k\right] }\left( J\left( 1,k\right) \right) \widehat{R}^{\left[ k%
\right] }\left( J\left( 1,k\right) \right) }{\rho _{A}\widehat{L}^{\left[ k%
\right] }\left( J\left( 0,k-1\right) \right) \widehat{R}^{\left[ k\right]
}\left( J\left( 0,k-1\right) \right) }=\frac{\nu \left( J\left( 0,k\right)
\right) }{\nu \left( J\left( 0,k-1\right) \right) }\mbox{.}
\end{eqnarray*}%
Hence $\nu $ is a $k$-step Markov measure.

\textbf{Step 3.} We finish the proof of Theorem \ref{Thm: 8} by
setting some notation first. For $k\geq 1$, let $A^{\left[ k\right]
}\in \mathbb{R}^{\mathbf{d}^k\times \mathbf{d}^k}$ be the $k$-th
higher block transition matrix which is indexed by $D^{[k]}$. Let
$\mathbb{A}=\mathbb{A}^{\left[ 1\right]
}=stoch\left( A\right) $ and $\mathbb{A}^{\left[ k\right] }=stoch\left( A^{%
\left[ k\right] }\right) $ for all $k>1$, and let $\mathbb{L}^{\left[ k\right] }$ and $\mathbf{1}_{a_{k}}$ be
the left and right eigenvectors of $\mathbb{A}^{\left[ k\right]}$ corresponding to the
eigenvalue $\rho _{\mathbb{A}^{\left[ k\right] }}=1$.

For $k\geq 1$, we define $S^{(0)}:=\left\langle \left\{ \mathbb{\widehat{L}}%
(i):i\in \mathcal{A}(Y)\right\} \right\rangle $ and
\begin{equation*}
S^{\left( k\right) } := \left\langle \left\{ \mathbb{\widehat{L}}^{\left[ k%
\right] }\left( J\left( 0,k-1\right) \right) \mathbb{N}_{J\left( 0,k\right)
}:J\left( 0,k\right) \in Y_{k+1}\right\} \right\rangle \mbox{,}
\end{equation*}%
where $\langle \;\; \rangle$ is used to denote span. Let
$\widehat{S}^{\left( k\right) }$ be the following sets.%
\begin{equation*}
\widehat{S}^{\left( 0\right) } = \left\langle \left\{ \left[ 0_{\mathbf{m}%
},\ldots ,0_{\mathbf{m}},\mathbb{\widehat{L}}(i),0_{\mathbf{m}},\ldots ,0_{%
\mathbf{m}}\right] :i\in \mathcal{A}(Y)\right\} \right\rangle
\end{equation*}%
and
\begin{equation*}
\widehat{S}^{\left( k\right) }=\left\langle \left\{
\begin{array}{c}
\left[ 0_{\mathbf{m}^{k}},\ldots
,0_{\mathbf{m}^{k}},\mathbb{\widehat{L}}^{\left[ k\right] }\left(
J\left( 0,k-1\right) \right) \mathbb{N}_{J\left( 0,k\right)
},0_{\mathbf{m}^{k}},\ldots ,0_{\mathbf{m}^{k}}\right] : \\
J\left( 0,k\right) \in Y_{k+1}%
\end{array}%
\right\} \right\rangle \mbox{.}
\end{equation*}%

We show  that if $\nu $ is an $n$-step Markov measure, then $%
n\leq \mathbf{m}-\mathbf{n}$. It is sufficient to prove that
\begin{equation}
S^{(k)} \subseteq S^{(k+1)} \hbox{ for all } k\geq 0.\label{33}
\end{equation}
This is true because if $\nu$ is an $n$-step Markov measure, then $S^{(n)} = S^{(n+1)}$.
Therefore,
\[
S^{(n)} = S^{(n+1)} \hbox{ implies } S^{(n)} = S^{(t)} \hbox{ for
all } t \geq n.
\]
Since $\dim(S^{(0)}) \leq \textbf{n}$, it follows that $n\leq
\textbf{m} - \textbf{n}$.

To prove (\ref{33}) we show that  $\widehat{S}^{(k)} \subseteq \widehat{S}^{(k+1)}$  for all  $k\geq 0$.

(\textbf{a}) For $k=0$, we claim $\widehat{S}%
^{(0)}\subseteq \widehat{S}^{(1)}$. Indeed,
\begin{equation*}
\widehat{S}^{\left( 1\right) } = \left\langle \left\{ %
\left[ 0_{\mathbf{m}},\ldots ,0_{\mathbf{m}},\mathbb{\widehat{L}}(i)\mathbb{%
N}_{ij},0_{\mathbf{m}},\ldots ,0_{\mathbf{m}}\right] :j\in \mathcal{A}%
(Y)\right\} \right\rangle \mbox{.}
\end{equation*}%
For all $j\in \mathcal{A}(Y)$, since $\mathbb{L}\mathbb{A}=\mathbb{L}$,
\begin{equation*}
\mathbb{\widehat{L}B}=\left( \mathbb{L}\left(P^{[k]}\right)^{-1}\right) \mathbb{B}P^{[k]}\left(P^{[k]}\right)^{-1}=%
\mathbb{L}\mathbb{A}\left(P^{[k]}\right)^{-1}=\mathbb{L}\left(P^{[k]}\right)^{-1}=\mathbb{\widehat{L}}.
\end{equation*}%
Thus,
\begin{eqnarray*}
\left[ 0_{\mathbf{m}},\ldots ,0_{\mathbf{m}},\mathbb{\widehat{L}}(j),0_{%
\mathbf{m}},\ldots ,0_{\mathbf{m}}\right] &=&\mathbb{\widehat{L}B}_{j}=(%
\mathbb{\widehat{L}B})\mathbb{B}_{j}=\sum_{i=1}^{\mathbf{n}}\mathbb{\widehat{%
L}B}_{i}\mathbb{B}_{j} \\
&=&\sum_{i=1}^{\mathbf{n}}\left[ 0_{\mathbf{m}},\ldots ,0_{\mathbf{m}},%
\mathbb{\widehat{L}}(i)\mathbb{N}_{ij},0_{\mathbf{m}},\ldots ,0_{\mathbf{m}}%
\right] .
\end{eqnarray*}%
Hence $\widehat{S}^{(0)}\subseteq \widehat{S}^{(1)}$.

(\textbf{b}) For $%
k\geq 1$, we claim $\widehat{S}^{(k)}\subseteq \widehat{S}^{(k+1)}$. Indeed, %
since $A^{[k]}$ and $A^{[k+1]}$ are the higher block transition matrices, it follows that $A^{[k]}$ and $A^{[k+1]}$ %
are shift equivalent (see \cite{lind-marcus-1995} for more detail), that is, there exists $F^{[k]}$ where a row is indexed by %
$D^{[k+1]}$ and a column is indexed by $D^{[k]}$ such that
\begin{equation}
F^{[k]}A^{[k]} = A^{[k+1]}F^{[k]}.\label{49}
\end{equation}
By recalling the definition of $stoch(A^{[k]})$, we have
\begin{equation}
\mathbb{A}^{[k]} = stoch(A^{[k]}) = \frac{1}{\rho_A} D_k^{-1} A^{[k]} D_k \label{43}
\end{equation}
where $D_k$ is the diagonal matrix with diagonal entries $D_k(i, i) = R^{[k]}(i), R^{[k]}$ is the right eigenvector of $A^{[k]}$ with spectral radius $\rho_{A^{[k]}} = \rho_A.$
Combining (\ref{49}) and (\ref{43}), yields
\begin{equation}
\mathbb{F}^{[k]}\mathbb{A}^{[k]} = \mathbb{A}^{[k+1]}\mathbb{F}^{[k]}\label{50}
\end{equation}
where $\mathbb{F}^{[k]} = D_{k+1}^{-1}F^{[k]}D_k$. It follows from
(\ref{50}) that we have
\begin{equation}
\mathbb{L}^{[k]}=\mathbb{L}^{[k+1]}\mathbb{F}^{[k]}\hbox{ for all }
k\geq 1\label{47}
\end{equation}
and
\begin{equation}
\mathbb{F}^{[k]}\mathbb{A}^{[k]}_{J(1, k)} = \left(\sum_{J(0,
k)}\mathbb{A}^{[k+1]}_{J(0, k)}\right)\mathbb{F}^{[k]} \hbox{ for
all } k\geq 1.\label{48}
\end{equation}
It follows from (\ref{47}) and (\ref{48}), for any $J(1,k+1)\in
Y_{k+1},$ we obtain
\begin{eqnarray}
&&\left[ 0_{a_{k}},\ldots ,0_{a_{k}},\mathbb{\widehat{L}}^{[k]}(J(1,k))%
\mathbb{N}_{J(1,k+1)},0_{a_{k}},\ldots ,0_{a_{k}}\right] \nonumber\\
&=&\left( \mathbb{L}^{[k+1]}\mathbb{F}^{[k]}\right) \mathbb{A}_{J(1,k)}^{[k]}\mathbb{A%
}_{J(2,k+1)}^{[k]}P^{[k]} \nonumber\\
&=&\mathbb{L}^{[k+1]}\left( \sum_{J(0,k)}\mathbb{A}_{J(0,k)}^{[k+1]}\right)
\left( \sum_{J(1,k+1)}\mathbb{A}_{J(1,k+1)}^{[k+1]}\right) \mathbb{F}^{[k]}P^{[k]} \nonumber\\
&=&\sum_{J(0,k)}\mathbb{L}^{[k+1]}\mathbb{A}_{J(0,k)}^{[k+1]}\mathbb{A}%
_{J(1,k+1)}^{[k+1]}P^{[k+1]}\label{51}
\end{eqnarray}%
Since $\mathbb{A}_J^{[k+1]}
=\left(P^{[k+1]}\right)^{-1}\mathbb{B}_J^{[k+1]} P^{[k+1]}$ for all
$k\geq 1, J\in Y_{k+1}$, we have
\begin{eqnarray}
&&\mathbb{L}^{[k+1]}\mathbb{A}_{J(0,k)}^{[k+1]}\mathbb{A}%
_{J(1,k+1)}^{[k+1]}P^{[k+1]} \nonumber\\
&=& \left[ 0_{\mathbf{m}^{k+1}},\ldots ,0_{\mathbf{m}^{k+1}},\widehat{\mathbb{L}}%
^{[k+1]}(J(0,k))\mathbb{N}_{J(0,k+1)},0_{\mathbf{m}^{k+1}},\ldots\label{52}
,0_{\mathbf{m}^{k+1}}\right]
\end{eqnarray}
On account of the above (\ref{51}) and (\ref{52}), we thus get
\begin{eqnarray*}
&&\left[ 0_{\mathbf{m}^k},\ldots ,0_{\mathbf{m}^k},\widehat{\mathbb{L}}^{[k]}(J(1,k))%
\mathbb{N}_{J(1,k+1)},0_{\mathbf{m}^k},\ldots ,0_{\mathbf{m}^k}\right] \\
&=&\sum_{J(0,k)}\left[ 0_{\mathbf{m}^{k+1}},\ldots ,0_{\mathbf{m}^{k+1}},\widehat{\mathbb{L}}%
^{[k+1]}(J(0,k))\mathbb{N}_{J(0,k+1)},0_{\mathbf{m}^{k+1}},\ldots
,0_{\mathbf{m}^{k+1}}\right].
\end{eqnarray*}%
Hence $\widehat{S}^{(k)}\subseteq \widehat{S}^{(k+1)}$ for all
$k\geq 1$.

According to (\textbf{a}) and (\textbf{b}),
we have $\widehat{S}^{(k)}  \subseteq \widehat{S}^{(k+1)}$ for all
$k\geq 0$. Hence, (\ref{33}) is proved, and so is the theorem.

\end{proof}

\subsection{Proof of Theorem \protect\ref{Thm: 3}}

To increase the readability we only prove the theorem for $k=1$, since the
general case is similar.

Since $X_M$ is SFT, we can construct the linear
representation of measure $\eta$ where $\eta$ is the unique maximal
measure on $X_M$. Combining the fact of $\widehat{L}=\left( \left[ \widehat{L}%
\left( k\right) \right] _{k=1}^{\mathbf{m}}\right) ^{t}$ is the left
eigenvector of $B$ corresponding to the maximal eigenvalue $\rho _{A}$ and $%
B^{\left( i,j\right) }=N_{ij}$ for all $i,j=1,\ldots ,\mathbf{n}$ yields
\begin{equation}
\sum_{i=1}^{\mathbf{n}}\widehat{L}\left( i\right) N_{ij}=\rho _{A}\widehat{L}%
\left( j\right) \mbox{ for all } j=1,\ldots ,\mathbf{n}.  \label{30}
\end{equation}%
Since $N$ satisfies the Markov condition from left of order $1$,
i.e., there exists $V_{1},\ldots ,V_{\mathbf{n}}$ such that
$V_{i}N_{ij}\subseteq V_{j}$ for $i,j=1,\ldots ,\mathbf{n}$. It can
be checked that $V_{k}=\widehat{L}\left( k\right) $ for $k=1,\ldots
,\mathbf{n}$ and $\widehat{L}\left( i\right)
N_{ij}=m\left( i,j\right) \widehat{L}\left( j\right) $ for $i,j=1,\ldots ,%
\mathbf{n}$. It follows from (\ref{30}), for all $j=1,\ldots ,\mathbf{n}$,
we have $\rho _{A}\widehat{L}\left( j\right) =\sum_{i=1}^{\mathbf{n}}%
\widehat{L}\left( i\right) N_{ij}=\sum_{i=1}^{\mathbf{n}}m\left( i,j\right)
\widehat{L}\left( j\right) $, which implies $\sum_{i=1}^{\mathbf{m}}m\left(
i,j\right) =\rho _{A}$ and $\mathbf{1}_{\mathbf{n}}^{t}M=\rho _{A}\mathbf{1}%
_{\mathbf{n}}^{t}$, where $\mathbf{1}_{\mathbf{n}}$ is $1\times \mathbf{n}$
column vector with entries which are all $1^{\prime }s$. This means that $\mathbf{1%
}_{\mathbf{n}}^{t}$ is a left eigenvector of $M$ corresponding to
eigenvalues $\rho_M = \rho _{A}$.

For any $k=1,\ldots ,\mathbf{n}$, we define $M_{k}$ which is the matrix indexed
by $\mathcal{A}\left( Y\right) =\left\{ 1,\ldots ,\mathbf{n}\right\}
$ as follows. For $i,j=1,\ldots ,\mathbf{n,}$
\begin{equation*}
M_{k}\left( i,j\right) =\left\{
\begin{array}{ccc}
m\left( i,j\right) &  & \mbox{if }j=k, \\
0 &  & \mbox{otherwise.}%
\end{array}%
\right.
\end{equation*}%
Choose $\widetilde{L}=\mathbf{1}_{\mathbf{n}}^{t}$ and
$\widetilde{R}=\left[
\widehat{L}\left( i\right) \widehat{R}\left( i\right) \right] _{i=1}^{%
\mathbf{n}}$, then the triple $\left(\widetilde{L},
\left\{M_i\right\}_{i=1}^{\textbf{n}}, \widetilde{R}\right)$
represents the unique maximal measure on $X_M$, that is, for any $J
= [j_0,\ldots, j_{n-1}] \in Y_n$,
\begin{equation}
\eta \left( J\right) :=\left( 1/\rho _{A}\right)
^{n}\widetilde{L}M_{j_{0}}\cdots M_{j_{n-1}}\widetilde{R}.\label{38}
\end{equation}

We claim that $\nu$ is the unique maximal measure
on $X_M$, i.e., $\nu(J) = \eta(J)$ for all $J \in Y_n, n
\in\mathbb{N}$. According to (\ref{38}), for any $J = [j_0,\ldots,
j_{n-1}] \in Y_n$, we get
\begin{eqnarray}
\eta \left( J\right) &=&\left( 1/\rho _{A}\right) ^{n}\widetilde{L}%
M_{j_{0}}\cdots M_{j_{n-1}}\widetilde{R} \;\;(by (\ref{38}))\nonumber\\
&=&\left( 1/\rho _{A}\right) ^{n-1}\rho _{A}\stackrel{1\mbox{ is at }j_{0}%
\mbox{-th coordinate}}{\overbrace{\left[ 0,\ldots ,0,1,0,\ldots ,0\right] }}%
M_{j_{1}}\cdots M_{j_{n-1}}\widetilde{R} \nonumber\\
&=&\left( 1/\rho _{A}\right) ^{n-1}m\left( j_{0},j_{1}\right) \left[
0,\ldots ,0,1,0,\ldots ,0\right] M_{j_{2}}\cdots
M_{j_{n-1}}\widetilde{R}.\label{36}
\end{eqnarray}%
Continue the same process as (\ref{36}), we have%
\begin{eqnarray*}
\eta \left( J\right) &=&\left( 1/\rho _{A}\right) ^{n-1}m\left(
j_{0},j_{1}\right) m\left( j_{1},j_{2}\right) \cdots m\left(
j_{n-2},j_{n-1}\right) \left[ 0,\ldots ,0,1,0,\ldots ,0\right] \widetilde{R}
\\
&=&\left( 1/\rho _{A}\right) ^{n-1}m\left( j_{0},j_{1}\right) m\left(
j_{1},j_{2}\right) \cdots m\left( j_{n-2},j_{n-1}\right) \widehat{L}\left(
j_{n-1}\right) \widehat{R}\left( j_{n-1}\right) .
\end{eqnarray*}%
For any $J\in Y_{n}$, we also have
\begin{eqnarray}
\nu \left( J\right) &=&\left( 1/\rho _{A}\right) ^{n}\rho _{A}\widehat{L}%
\left( j_{0}\right) N_{j_{0}j_{1}}N_{j_{1}j_{2}}\cdots N_{j_{n-2}j_{n-1}}%
\widehat{R}\left( j_{n-1}\right) \nonumber\\
&=&\left( 1/\rho _{A}\right) ^{n}m\left( j_{0},j_{1}\right) \widehat{L}%
\left( j_{1}\right) N_{j_{1}j_{2}}\cdots
N_{j_{n-2}j_{n-1}}\widehat{R}\left( j_{n-1}\right) .\label{37}
\end{eqnarray}%
Continue the same process as (\ref{37}), we have
\begin{eqnarray*}
\nu \left( J\right) &=&\left( 1/\rho _{A}\right) ^{n}m\left(
j_{0},j_{1}\right) m\left( j_{1},j_{2}\right) \cdots m\left(
j_{n-2},j_{n-1}\right) \widehat{L}\left( j_{n-1}\right) \widehat{R}\left(
j_{n-1}\right) \\
&=&\eta \left( J\right) \mbox{.}
\end{eqnarray*}%
Hence, $\nu $ is the unique maximal measure of $X_M$. This completes
the proof.

\subsection{Proof of Theorem \protect\ref{Thm: 6}}

We present some results from \cite{Kenyon1996} for the criterion of the
equality for Hausdorff and Minkowski dimensions.

\begin{theorem}[{\cite[Theorem 1.3]{Kenyon1996}}]
\label{Thm: 13}Under the same assumptions of Theorem \ref{Thm: 1} with $A$
being irreducible. Then $\dim _{H}Z=\dim _{M}Z$ iff the unique invariant
measure of maximal entropy in $Z$ projects via $\pi _{y}$ to the unique
measure of maximal entropy on $\pi _{y}(Z)$.
\end{theorem}

\begin{proof}[Proof of (i) of Theorem \protect\ref{Thm: 6}]
($\Rightarrow $) Combining the fact that $N$ satisfies the Markov
condition of order $k$ and Theorem \ref{Thm: 3}, we have $\pi \mu
_{A}=\nu $ is the unique maximal measure of subshift of finite type
$X_M$ with adjacency matrix $M$ where $\mu_A$ is the measure of
maximal entropy in $Z$. If $\rho_{M}=\rho_{T^{[k]}}$, then $\nu$ is
also the unique maximal measure of $Y$. Hence $\dim _{H}Z=\dim
_{M}Z$ by Theorem \ref{Thm: 13}.

($\Leftarrow $) If $\dim _{H}Z=\dim _{M}Z$, then the unique invariant
measure of maximal entropy in $Z$ projects via $\pi _{y}$ to the unique
measure of maximal entropy on $\pi _{y}(Z)$. Hence, $\rho_{M}=\rho_{T^{[k]}}$.
\end{proof}

\begin{proof}[Proof of (ii) of Theorem \protect\ref{Thm: 6}]
Without loss of generality, we may assume $N$ satisfies the Markov condition of
order 1 from left. Let $A^{\prime }$ be a transition matrix of $Z^{\prime
}=K(\mathbb{T},D^{\prime })$. For all $k=1,\ldots,\mathbf{n},$ we denote $%
A^{\prime }_k$ by
\begin{equation*}
A^{\prime }_k\left( d,d^{\prime }\right) =\left\{
\begin{array}{ccc}
1, &  & \mbox{if }d^{\prime }=(d^{\prime }_1,d^{\prime }_2) \mbox{ with }%
d^{\prime }_2=k;\\
0, &  & \mbox{otherwise,}%
\end{array}%
\right.
\end{equation*}
where $d, d^{\prime }\in D^{\prime }$.
Then
\begin{eqnarray*}
&&\dim _{H}Z^{\prime }=\dim _{H}K(\mathbb{T},D^{\prime }) \\
&=&\lim_{N\rightarrow \infty }\frac{1}{N}\log _{\mathbf{n}}\sum_{J(1,N)\in
Y_{N}}\left\Vert A_{j_{1}}^{\prime }A_{j_{2}}^{\prime }\cdots
A_{j_{N}}^{\prime }\right\Vert ^{\alpha } \\
&=&\lim_{N\rightarrow \infty }\frac{1}{N}\log _{\mathbf{n}}\sum_{J(1,N)\in
Y_{N}}\left( LA_{j_{1}}A_{j_{2}}\cdots A_{j_{N}}R\right) ^{\alpha } \\
&=&\lim_{N\rightarrow \infty }\frac{1}{N}\log _{\mathbf{n}}\sum_{J(1,N)\in
Y_{N}}\left( \rho _{A}m(j_{1},j_{2})\cdots m(j_{N-1},j_{N})\widehat{L}(j_{N})%
\widehat{R}(j_{N})\right) ^{\alpha } \\
&=&\lim_{N\rightarrow \infty }\frac{1}{N}\log _{\mathbf{n}}\left( \rho
_{A}^{\alpha }\mathbf{1}_{\mathbf{n}}^{t}\left( M^{\alpha }\right) ^{N-1}%
\left[ \widehat{L}(i)\widehat{R}(i)\right] _{i=1}^{\mathbf{n}}\right) \\
&=&\lim_{N\rightarrow \infty }\frac{1}{N}\log _{\mathbf{n}}\left\Vert \left(
M^{\alpha }\right) ^{N-1}\right\Vert
\end{eqnarray*}%
From the proof of Theorem \ref{Thm: 8}, we know $N_{ij}=\left[%
A\left(\left(p, i\right), \left(q, j\right)\right)\right]_{p, q=1}^{\mathbf{m}%
}$ and $V_i = \widehat{L}(i) = \left(\left[L(d_1, i)\right]_{d_1=1}^{\mathbf{m}%
}\right)^t$ for all $i, j =1,\cdots,\mathbf{n}$.  Recall that $z(i)$ is the number of rectangles in row $i$. Since
\begin{equation}
\label{54}
L(d_1, i) =\left\{
\begin{array}{ccc}
1, &  & \mbox{if }(d_1, i)\in D^{\prime }\mbox{ } ;\\
0, &  & \mbox{otherwise}%
\end{array}%
\right.
\end{equation}
and
\begin{equation}
\label{55}
A\left(\left(p, i\right), \left(q, j\right)\right) =\left\{
\begin{array}{ccc}
1, &  & \mbox{if }\left(p, i\right), \left(q, j\right)\in D^{\prime
};
\\
0, &  & \mbox{otherwise},%
\end{array}%
\right.
\end{equation}
we obtain for all $i = 1,\cdots,\mathbf{n}$
\begin{equation}
\label{53}\sum_{(p, i)\in D'} A\left((p, i), (q, j)\right) = z(i).
\end{equation}
According to (\ref{54}), (\ref{55}) and  (\ref{53}), we thus get for all $i,j = 1,\cdots,\mathbf{n}$
\begin{eqnarray*}
&&V_iN_{ij} = \widehat{L}(i)N_{i,j} \\
&=&\left[\sum_{p=1}^{\mathbf{m}}L(p,i)A\left(\left(p,i\right),\left(q,j%
\right)\right)\right]_{q=1}^{\mathbf{m}}=\left[\sum_{(p,i)\in D^{\prime
}}A\left(\left(p,i\right),\left(q,j\right)\right)\right]_{q=1}^{\mathbf{m}}
\\
&=&\sum_{(p,i)\in D^{\prime
}}A\left(\left(p,i\right),\left(q,j\right)\right)\left[L(q,j)\right]_{q=1}^{%
\mathbf{m}}=z(i)\left[L(q,j)\right]_{q=1}^{\mathbf{m}}=m(i,j)V_j,
\end{eqnarray*}
that is, $m(i,j)=z(i)$ for all $i,j=1,\ldots ,\mathbf{n}$. Thus, $\log _{%
\mathbf{n}}\left\Vert \left( M^{\alpha }\right) ^{N-1}\right\Vert =\log _{%
\mathbf{n}}\sum_{j=1}^{\mathbf{n}}z(j)^{\alpha }$. The proof is thus
completed.
\end{proof}

\begin{proof}[Proof of (iii) of Theorem \protect\ref{Thm: 6}]
For $k\geq 1$, it follows from the same argument of $k=1$, that we have
\begin{equation*}
\dim _{H}\left[ K_{T}\left( A\right) \right] =\lim_{N\rightarrow \infty }%
\frac{1}{N}\log _{\mathbf{n}}\sum_{J\left( 0,N-1\right) \in Y_{N}}\left[ L^{%
\left[ k\right] }\left( \prod_{i=0}^{N-k}A_{J\left( i,i+k-1\right) }^{\left[
k\right] }\right) R^{\left[ k\right] }\right] ^{\alpha }\mbox{.}
\end{equation*}%
And for any $J=J\left( 0,N-1\right) \in Y_{N}$ with $N\geq k$, we have
\begin{eqnarray*}
&&L^{\left[ k\right] }\left( \prod_{i=0}^{N-k}A_{J\left( i,i+k-1\right) }^{%
\left[ k\right] }\right) R^{\left[ k\right] }=\rho _{A}\left(
\prod_{i=0}^{N-k-1}m\left( J\left( i,i+k-1\right) ,J\left( i+1,i+k\right)
\right) \right) \times \\
&&\widehat{L}^{\left[ k\right] }\left( J\left( N-k-1,N-1\right) \right)
\widehat{R}^{\left[ k\right] }\left( J\left( N-k-1,N-1\right) \right) .
\end{eqnarray*}%
Hence,
\begin{eqnarray*}
&&\dim _{H}Z \\
&=&\lim_{N\rightarrow \infty }\frac{1}{N}\log _{\mathbf{n}}\sum_{J\left(
0,N-1\right) \in Y_{N}}[\rho _{A}\left( \prod_{i=0}^{N-k-1}m\left( J\left(
i,i+k-1\right) ,J\left( i+1,i+k\right) \right) \right) \times \\
&&\widehat{L}^{\left[ k\right] }\left( J\left( N-k-1,N-1\right) \right)
\widehat{R}^{\left[ k\right] }\left( J\left( N-k-1,N-1\right) \right)
]^{\alpha } \\
&=&\lim_{N\rightarrow \infty }\frac{1}{N}\log
_{\mathbf{n}}\left(\rho_A^{\alpha}\textbf{1}_{a_k}^t\left(M^{\alpha}\right)^{N-k-2}\left[\widehat{L}^{[k]}(J_p)\widehat{R}^{[k]}(J_p)\right]_{p=1}^{a_k}\right)\\
&=&\lim_{N\rightarrow \infty }\frac{1}{N}\log
_{\mathbf{n}}\left\Vert \left(
M^{\alpha }\right) ^{N-k-2}\right\Vert =\log _{\mathbf{n}}\rho _{M^{\alpha }}%
\mbox{.}
\end{eqnarray*}%
The proof is completed.
\end{proof}
\section{Application and examples}

\subsection{Criterion for Markov measure}

Assume $\pi :X\rightarrow Y$ is an one-block code induced from $\Pi :%
\mathcal{A}(X)\rightarrow \mathcal{A}(Y)$, we assume that $\mathcal{A}%
(Y)=\left\{ 1,\ldots ,\mathbf{n}\right\} $ and $X$ is a SFT with the
transition matrix $A$ which is irreducible. Suppose $Y$ is a
irreducible subshift of finite type, we call $\pi :X\rightarrow Y$ \emph{Markovian}
(cf.~\cite{Boyle2011, BoyleTuncel1984}) if for every Markov measure $\nu $ on $Y$%
, there is a Markov measure $\mu $ on $X$ with $\pi \mu =\nu $.

For $k=1,\ldots ,\mathbf{n}$, let%
\begin{equation*}
E_{k}=\left\{ i\in \mathcal{A}(X):\Pi \left( i\right) =k\right\}
\mbox{,}
\end{equation*}%
and denote by $e_{k}=\#E_{k}$ the number of $E_{k}$ and define
$\widehat{N}_{ij} \in
\mathbb{R}^{e_{i}\times e_{j}}$ as
follows:%
\begin{equation*}
\widehat{N}_{ij}(k,l)=\left\{
\begin{array}{ccc}
1, &  & \mbox{if }A(k,l)=1\mbox{ and }\Pi (k)=i\mbox{ and }\Pi (l)=j ;\\
0, &  & \mbox{otherwise}.%
\end{array}%
\right.
\end{equation*}%
For $1\leq i,j\leq \mathbf{n}$, $\widehat{N}_{ij}$ is called
\emph{row allowable }if for each $k\in \left\{ 1,\ldots
,e_{i}\right\} $ there is an $l\in \left\{ 1,\ldots ,e_{j}\right\} $
such that $\widehat{N}_{ij}(k,l)\neq 0$. $\widehat{N}=\left(
\widehat{N}_{ij}\right) _{i,j=1}^{\mathbf{n}}$ is called \emph{row allowable }if $%
\widehat{N}_{ij}$ is row allowable for all $1\leq i,j\leq \mathbf{n}$. In \cite%
{Chazottes2003}, the authors call such factors \emph{full row
allowable}. The following lemma shows that the full row
allowability implies the projection space $Y$ is a subshift of
finite type.

\begin{lemma}[{\cite[Lemma 6]{Chazottes2003}}]
\label{Lma: 4}If $\pi $ is full row allowable, then $Y$ is a
subshift of finite type.
\end{lemma}

Define (\textbf{H1}) and (\textbf{H2}) as follows:

\begin{enumerate}
\item[(\textbf{H1})] $\widehat{N}=\left( \widehat{N}_{ij}\right) _{i,j=1}^{\mathbf{n}}$ is row
allowable.

\item[(\textbf{H2})] For all $J\in Per_{n}\left( Y\right) $, the $n$%
-periodic orbit in $Y$, with $1\leq n\leq \#\mathcal{A}(Y)$,
$\widehat{N}_{J}$ is a positive matrix.
\end{enumerate}

Under (\textbf{H1}) and (\textbf{H2}), Chazottes and Ugalde prove
that there is a Gibbs measure of some well-defined potential on $Y$.
In the following, we give another proof for this result and give a
simple lemma first.

\begin{lemma}
\label{Lma: 3}Let $E_{n}\in R^{n\times n}$ be the full matrix and $%
C_{1},C_{2},\ldots ,C_{m}$ be a sequence of row allowable matrices, then $%
C_{1}C_{2}\cdots C_{m}E_{n}\geq E_{n}$.
\end{lemma}

\begin{proof}
This is the immediate consequence of the observation that
$AE_{n}\geq E_{n}$ if $A$ is row allowable.
\end{proof}

\begin{theorem}
\label{Thm: 11} Let $\pi :X\rightarrow Y$ be an one-block code which satisfies (%
\textbf{H1}) and (\textbf{H2}), then there exists a unique Gibbs measure $%
\nu $ on $Y$.
\end{theorem}

\begin{proof}
Without loss of generality, we may assume that $\widehat{N}=\left(
\widehat{N}_{ij}\right) _{i,j=1}^{\mathbf{n}}$ are all square
matrices with the same size $\mathbf{m}
$, that is, $\widehat{N}_{ij}\in \mathbb{R}^{\mathbf{m}\times \mathbf{m}}$ for all $%
1\leq i,j\leq \mathbf{n}$. We first claim that if $\pi $ satisfies (%
\textbf{H1}) and (\textbf{H2}), then $f\in \mathcal{D}_{w}(Y,p)$ for some $%
p\in \mathbb{N}$, where $f(I)=\left\Vert \widehat{N}_{I}\right\Vert
$ for all $I\in Y_{n}$. Indeed, since $\left\Vert
\widehat{N}_{IJ}\right\Vert \leq \left\Vert
\widehat{N}_{I}\right\Vert \left\Vert \widehat{N}_{J}\right\Vert $,
we only need to prove that there is a $p\in \mathbb{N}$ such that
for all $I$, $J\in Y^{\ast }$, there exists $K\in \cup
_{i=0}^{p}Y_{i}$ with $IKJ\in Y^{\ast }$ and $\left\Vert
\widehat{N}_{IKJ}\right\Vert \geq c\left\Vert
\widehat{N}_{I}\right\Vert \left\Vert \widehat{N}_{J}\right\Vert $
for some $c>0$. For $I=\left( i_{0},\ldots ,i_{m-1}\right) \in
Y_{m}$ and $J=\left( j_{0},\ldots ,j_{n-1}\right) \in Y_{n}$, take
$\bar{i}_{m-1}$, $\bar{j}_{0}\in \mathcal{A}(X)$ such that
$\Pi (\bar{i}_{m-1})=i_{m-1}$ and $\Pi
(\bar{j}_{0})=j_{0}$. Since $X$ is irreducible, there is a path
in $X$ of the form $\bar{i}_{m-1}$
to $\bar{j}_{0}$ of length $m_{1}$, say $\bar{I}^{\prime }=\left( \bar{i}_{m-1},\ldots ,\bar{j}_{0}\right) $, and a periodic path from
$\bar{j}_{0}$
to $\bar{j}_{0}$ of length $m_{2}$, say $\bar{P}=\left( \bar{j}_{0},\ldots ,%
\bar{j}_{0}\right) $. Denote by $\bar{I}^{\prime }\bar{P}=\left( \bar{i}%
_{m-1},\ldots ,\bar{j}_{0},\ldots ,\bar{j}_{0}\right) $ the
concatenation of $\bar{I}$ and $\bar{P}$ and let $I^{\prime }=\pi
(\bar{I}^{\prime })=\Pi \left( \bar{i}_{m-1}\right) \cdots \Pi
\left( \bar{j}_{0}\right) \in Y_{m_{1}}$ and $P=\pi (\bar{P})=\Pi
\left( \bar{j}_{0}\right) \cdots \Pi \left( \bar{j}_{0}\right) \in
Y_{m_{2}}$. Since $\pi (\bar{P})$ is a
periodic path in $Y_{m_{2}}$, (\textbf{H1}) is applied to show that $%
\widehat{N}_{P}\geq E_{\mathbf{m}}$, then Lemma \ref{Lma: 3} and
(\textbf{H2}) is applied to obtain
\begin{equation*}
\left\Vert \widehat{N}_{II^{\prime }PJ}\right\Vert  =\left\Vert
\widehat{N}_{I}\widehat{N}_{I^{\prime
}}\widehat{N}_{P}\widehat{N}_{J}\right\Vert \geq \left\Vert
\widehat{N}_{I}E_{\mathbf{m}}\widehat{N}_{J}\right\Vert \geq
c\left\Vert \widehat{N}_{I}\right\Vert \left\Vert
\widehat{N}_{J}\right\Vert .
\end{equation*}%
Since both the length of $I^{\prime }$ and $P$ can be chosen so as to be less than $\mathbf{n}$, i.e., the number of $\mathcal{A}(Y)$. Then the claim follows if we take $%
p=2\mathbf{m}$. Hence, it follows from Theorem 5.5 of
\cite{Feng2011} that there is a unique equilibrium measure $\nu $
which is ergodic and satisfies
the following Gibbs property, that is, there is a $c >0$ such that for all $%
n\in \mathbb{N}$ and $J\in Y_{n}$
\begin{equation*}
c^{-1}\leq \frac{\nu (J)}{\exp \left( -nP(\sigma
_{Y},\widehat{N})\right) f(J)}\leq c.
\end{equation*}%
The proof is thus completed.
\end{proof}

A further question arose in \cite{Chazottes2003}: \emph{When is the
factor map not a Markov map? }Recall the result of Boyle and Tuncel
for the criterion of the Markovian factor $\pi $.

\begin{theorem}[\protect\cite{BoyleTuncel1984}]
\label{Thm: 15}For a factor map $\pi :X\rightarrow Y$ between
irreducible SFTs, if there exists any fully supported Markov measure
$\mu $ and $\nu $ with $\pi \mu =\nu ,$ then $\pi $ is Markovian.
\end{theorem}

We use the skill in Theorem \ref{Thm: 8} to answer the above
question up to minor modification of the induced potential $N$.
Arrange the set $E_{k}$ as an ordered set $\left\{ i_{j}^{k}:\Pi
\left( i_{j}^{k}\right) =k\right\} _{j=1}^{e_{k}}$, define $\left(
j,k\right) =i_{j}^{k}$ for all $1\leq k\leq \mathbf{n}$ and $1\leq
j\leq e_{k}$. Let $\mathbf{m}=\max_{1\leq k\leq \mathbf{n}}\left\{
e_{k}\right\} $. Introduce new symbols $D=\left\{ 1,\ldots
,\mathbf{m}\right\} \times \left\{ 1,\ldots ,\mathbf{n}\right\} $
and denote $\mathbf{d}=\mathbf{m}\times \mathbf{n}$, define the \emph{%
modified transition matrix} $B\in \mathbb{R}^{\mathbf{d}\times \mathbf{%
d}}$ which is indexed by $D$ as follows.%
\begin{equation*}
B(d,d^{\prime })=\left\{
\begin{array}{ccc}
1, &  & \mbox{if }d=\left( j,k\right) \in E_{k}\mbox{, }d^{\prime
}=\left(
j^{\prime },k^{\prime }\right) \in E_{k^{\prime }}\mbox{ with }%
A(i_{j}^{k},i_{j^{\prime }}^{k^{\prime }})=1 ;\\
0, &  & \mbox{otherwise}.%
\end{array}%
\right.
\end{equation*}%
Let $N$ be the induced matrix-valued potential from $B$, and call
$N$ the \emph{modified induced (matrix-valued) potential} on $Y $.
Note here that $N$ is not full row allowable, however, $\widehat{N} = \left(
N_{ij}(u,v)|_{u\in E_{i},v\in E_{j}}\right) _{i,j=1}^{\mathbf{n}}$
is row allowable for each $1\leq i,j\leq \mathbf{n}$ from
(\textbf{H1}). The following result comes from \cite{Heller1965}
which provides a criterion for Markov measures by means of a reduced
 module. To avoid the notation abuse, we omit the
definitions of reduced  module for measures and refer to
\cite{BarralFeng2009, Heller1965} and some references
therein.

\begin{theorem}[{\cite[Theorem 5.1]{Boyle2011}}, {\cite[Proposition 3.2]{Heller1965}}]
\label{Thm: 16}Let $\left( l,M,r\right) $ be a presentation of the
reduced  module of a sofic measure $\nu $ on $Y$, in which $M_{i}$
denotes the matrix by which a symbol $i$ of $\mathcal{A}\left(
Y\right) $ acts on the module. Suppose $k\in \mathbb{N}.$ Then the
sofic measure $\nu $ is $k$-step Markov if and only if every product
$M_{i\left( 1\right) }\cdots M_{i\left( k\right) }$ of length $k$
has a rank of $1$ at most.
\end{theorem}

\begin{theorem}\label{Thm: 12}
If $\widehat{N}$ satisfies (\textbf{H1%
}) and (\textbf{H2}), then the projection measure $\nu $ is Markov if
and only if $N$ satisfies the Markov condition where $N=\left( N_{ij}\right) _{i,j=1}^{\mathbf{n}}$ is
the modified induced matrix-valued potential on $Y$.
\end{theorem}

\begin{proof}
According to Theorem \ref{Thm: 8}, it is sufficient to show that if the projection measure $\nu $ is Markov, then $N$ satisfies the Markov condition.  Without loss of generality, we may assume that $\nu
$ is a 1-step Markov measure. According to Proposition \ref{Prop: 3} (ii), we have a natural presentation of a module $\left(L, \{B_{i}\}_{i=1}^{\mathbf{n}}, R\right)$ of a projection measure $\nu $ on $Y$, that is, let $B\in \mathbb{R}^{\mathbf{d}\times \mathbf{d}} \left(\mathbf{d}=\mathbf{m} \times%
\mathbf{n}\right)$ be the modified transition matrix and $L, R$ be
the left and right eigenvector of $B$ corresponding to the maximal
eigenvalue $\rho_{B}.$ Then
$\left(L,\{B_{i}\}_{i=1}^{\mathbf{n}},R\right)$ is a presentation of
a module of a projection measure $\nu $ on $Y$, where $B_{i}$ is
defined by
\begin{equation*}
B_{i}\left( d,d^{\prime }\right) =\left\{
\begin{array}{ccc}
B(d,d^{\prime }), &  & \mbox{if }d^{\prime }=(d_{1}^{\prime },d_{2}^{\prime })%
\mbox{ with }d_{2}^{\prime }=i; \\
0, &  & \mbox{otherwise},%
\end{array}%
\right.
\end{equation*}%
for all $d,d^{\prime }\in D$ and $i=1,\ldots ,\mathbf{n}$.

If $rank(B_i)=1$ for all $i\in\mathcal{A}(Y)$, by the same argument of the proof in Theorem \ref{Thm: 8}, we have $N$ which satisfies the Markov condition. If not, we can construct a smaller module $(l,\{\widetilde{B}_{i}\}_{i=1}^{\mathbf{n}},r)$ via $(L,\{B_{i}\}_{i=1}^{\mathbf{n}},R)$ such that $(l,\{\widetilde{B}_{i}\}_{i=1}^{\mathbf{n}},r)$ is the reduced module of $\nu$ on $Y$. This follows the same method as in \cite{Boyle2011}.

Indeed, let $\mathcal{U}$ be the vector space generated by vectors of the form $$
LB_{J(0,n-1)}=LB_{j_{0}}B_{j_{1}}\cdots B_{j_{n-1}},$$ for all $
J(0,n-1)\in Y_{n},n\in \mathbb{N}.$ Let $\dim \mathcal{U} = k$. If $k < \textbf{d}$, then construct a smaller module (presenting the same measure) as follows. Let
\begin{equation*}
\mathcal{B}:=\left\{ u_{i}=LB_{J(0,n-1)}\mbox{ for some }J(0,n-1)\in
Y_{n},n\in \mathbb{N}:i=1,\ldots ,k\right\}
\end{equation*}
be a basis of $\mathcal{U}$. By Lemma \ref{Lma: 6}, for each $a\in \mathcal{A}(Y)$, we
define ordered sets
\begin{equation*}
K_{a}:=\left\{
\begin{array}{c}
u_{l}^{(a)}=\left[ 0_{\mathbf{m}},\ldots ,0_{\mathbf{m}},v_{l}^{(a)},0_{%
\mathbf{m}},\ldots ,0_{\mathbf{m}}\right] \in \mathbb{R}^{1\times \mathbf{d}%
},v_{l}^{(a)}\in \mathbb{R}^{1\times \mathbf{m}}: \\
u_{l}^{(a)}=LB_{J(0,n-1)}\in \mathcal{B}\mbox{ and }j_{n-1}=a%
\end{array}%
\right\}
\end{equation*}%
and
\begin{equation*}
K_{a}^{\prime }:=\left\{ v_{l}^{(a)}:l=1,\ldots ,k_{a}\right\}
\mbox{, where }k_{a}=\#K_{a}.
\end{equation*}%
The matrix $U=(U^{(i,j)})_{i,j=1}^{\mathbf{n}}$ arising from $K_{a}$
is
defined by $U^{(i,j)}=0_{k_{i}\times \mathbf{m}}$ if $i\neq j$ and $%
U^{(i,j)} $ whose rows form an ordered set $K_{i}^{\prime }$ if
$i=j$. Then $UB_{i}=\widetilde{B}_{i}U$ for all $i=1,\ldots ,%
\mathbf{n}$, and $lU=L, r=UR$. Thus, $\rho _{B}=\rho
_{\widetilde{B}}$, where $\widetilde{B} = \sum_{i=1}^{\mathbf{n}} \widetilde{B}_i$. According to the form of $U$, we obtain
\begin{eqnarray}
& &\left[ 0_{k_{1}},\ldots ,0_{k_{j-1}},l(j),0_{k_{j+1}},\ldots
,0_{k_{\mathbf{n}}}\right] U\nonumber\\
&=&\left[ 0_{\mathbf{m}},\ldots ,0_{\mathbf{m}%
},(lU)(j),0_{\mathbf{m}},\ldots ,0_{\mathbf{m}}\right]\nonumber\\
&=&\left[ 0_{\mathbf{m}},\ldots ,0_{\mathbf{m}%
},L(j),0_{\mathbf{m}},\ldots ,0_{\mathbf{m}}\right]\label{45}
\end{eqnarray}

Without loss of generality, we may assume that $(l,\{\widetilde{B}_{i}\}_{i=1}^{\mathbf{n}},r)$ is the reduced module of $\nu$ on $Y$.  By Theorem \ref{Thm: 16}, we
obtain rank$(\widetilde{B}_{i})=1$ for all $i=1,\ldots ,\mathbf{n}.$ Let $%
l=[l(1),\ldots ,l(\mathbf{n})]$ where $l(i)\in \mathbb{R}^{1\times
k_{i}}$ for all $i$. Under the same argument of the proof in Theorem
\ref{Thm: 8}, the induced matrix-valued potential $\widetilde{N}$
satisfies the Markov condition, that is,
\begin{equation}
l(i)\widetilde{N}_{ij}=\widetilde{m}(i,j)l(j), \label{46}
\end{equation}
where $%
\widetilde{N}_{ij}\in \mathbb{R}^{k_{i}\times k_{j}}$ and $\widetilde{m}%
(i,j)\in \mathbb{R}$. Let $L=\left( \left[ L(i)\right] _{i=1}^{\mathbf{n}%
}\right) ^{t}$ where $L(i)=\left( [L(d_{1},i)]_{d_{1}\in D}\right)
^{t}$. Combining (\ref{45}),  (\ref{46}), for any $[ij]\in Y_{2},$
we have
\begin{eqnarray*}
&&LB_{i}B_{j}=\left( lU\right) B_{i}B_{j}=l\widetilde{B}_{i}\widetilde{B}%
_{j}U \\
&=&\left[ 0_{k_{1}},\ldots ,0_{k_{j-1}},\rho _{\widetilde{B}}l(i)\widetilde{N%
}_{ij},0_{k_{j+1}},\ldots ,0_{k_{\mathbf{n}}}\right] U \\
&=&\rho _{\widetilde{B}}\left[ 0_{k_{1}},\ldots ,0_{k_{j-1}},\widetilde{m}%
(i,j)l(j),0_{k_{j+1}},\ldots ,0_{k_{\mathbf{n}}}\right] U \\
&=&\rho _{B}\widetilde{m}(i,j)\left[ 0_{\mathbf{m}},\ldots ,0_{\mathbf{m}%
},L(j),0_{\mathbf{m}},\ldots ,0_{\mathbf{m}}\right]
\end{eqnarray*}%
On the other hand, since $LB_{i}B_{j}=\rho _{B}\left[ 0_{\mathbf{m}%
},\ldots ,0_{\mathbf{m}},L(i)N_{ij},0_{\mathbf{m}},\ldots ,0_{\mathbf{m}}%
\right] $, then we have $L(i)N_{ij}=\widetilde{m}(i,j)L(j)$ for all
$[ij]\in Y_{2}$. Hence we pick $V_i = L(i)$ for all $i = 1, \ldots,
\textbf{n}$ and the proof is thus completed.
\end{proof}
\subsection{Examples}

We give two examples illustrating Theorem \ref{Thm: 8}, Theorem
\ref{Thm: 3} and the application of the Markovian property for a
factor $\pi $.

In view of Theorem \ref{Thm: 16}, readers may wonder whether $\left( L,\left\{ B_{i}\right\}_{i=1}^{\textbf{n}} ,%
R\right) $ is one reduced module of $Y$. The following example
demonstrates that $N$ satisfies the Markov condition, however,
$\left(L, \left\{ B_{i}\right\}_{i=1}^{\textbf{n}}, R\right) $ is
not a reduced module.

\begin{example}[Blackwell]
Let $\mathcal{A}(X)=\left\{ 1,2,3\right\} $ and
$\mathcal{A}(Y)=\left\{
1,2\right\} $ and the one-block map $\Pi :\mathcal{A}(X)\rightarrow \mathcal{%
A}(Y)$ be defined by $\Pi (1)=1,$ $\Pi \left( 2\right) =2$ and $\Pi
\left( 3\right) =2$. Let $\pi :X\rightarrow Y$ be the factor from
$X$ to $Y$ induced from $\Pi $ with $X=\Sigma _{A}$ for some
\begin{equation*}
A=\left(
\begin{array}{ccc}
0 & 1 & 1 \\
1 & 1 & 0 \\
1 & 0 & 1%
\end{array}%
\right) \mbox{.}
\end{equation*}%
This factor has been proven ({\cite[Example 2.7]{Boyle2011}}) to be
Markovian. Here we use Theorem \ref{Thm: 12} to give a criterion for
this property. Since
$E_{1}=\left\{ 1\right\} $ and $E_{2}=\left\{ 2,3\right\} $ we see that $%
\mathbf{m}=e_{2}=2$, and we introduce the new symbols and the
corresponding sets $\widehat{E}_{1}$ and $\widehat{E}_{2}$ are as
follows.
\begin{eqnarray*}
D &=&\left\{ 1,2\right\} \times \left\{ 1,2\right\} =\left\{ \left(
1,1\right) ,\left( 1,2\right) ,\left( 2,1\right) ,\left( 2,2\right)
\right\}
, \\
\mbox{ }\widehat{E}_{1} &=&\left\{ 1=\left( 1,1\right) ,\left(
2,1\right) \right\} \mbox{, }\widehat{E}_{2}=\left\{ 2=\left(
1,2\right) ,3=\left( 2,2\right) \right\} .
\end{eqnarray*}%
Therefore,
\begin{equation*}
B = \left(
\begin{array}{cc}
N_{11} & N_{12} \\
N_{21} & N_{22}%
\end{array}%
\right) =\left(
\begin{array}{cccc}
0 & 0 & 1 & 1 \\
0 & 0 & 0 & 0 \\
1 & 0 & 1 & 0 \\
1 & 0 & 0 & 1%
\end{array}%
\right) .
\end{equation*}%
\begin{equation*}
N_{11}=\left(
\begin{array}{cc}
0 & 0 \\
0 & 0%
\end{array}%
\right) \mbox{, }N_{12}=\left(
\begin{array}{cc}
1 & 1 \\
0 & 0%
\end{array}%
\right) \mbox{, }N_{21}=\left(
\begin{array}{cc}
1 & 0 \\
1 & 0%
\end{array}%
\right) \mbox{, }N_{22}=\left(
\begin{array}{cc}
1 & 0 \\
0 & 1%
\end{array}%
\right) \mbox{.}
\end{equation*}%
Take $V_{1}=\left(
\begin{array}{cc}
1 & 0%
\end{array}%
\right) $ and $V_{2}=\left(
\begin{array}{cc}
1 & 1%
\end{array}%
\right) $, one can easily check that $N=\left( N_{ij}\right)
_{i,j=1}^{2}$ satisfies the Markov condition of order $1$. Thus
Theorem \ref{Thm: 12} is applied to show that the factor is a Markov
map. However, we remark here that $\left(L, \left\{
B_{i}\right\}_{i=1}^{\textbf{n}}, R\right) $ is not a linear
representation of the reduced module for sofic measure $\nu $ on $%
Y $ since rank$\left( N_{22}\right) =2$. Let
\begin{equation*}
B_{1}=\left(
\begin{array}{cccc}
0 & 0 & 0 & 0 \\
0 & 0 & 0 & 0 \\
1 & 0 & 0 & 0 \\
1 & 0 & 0 & 0%
\end{array}%
\right) \mbox{ and }B_{2}=\left(
\begin{array}{cccc}
0 & 0 & 1 & 1 \\
0 & 0 & 0 & 0 \\
0 & 0 & 1 & 0 \\
0 & 0 & 0 & 1%
\end{array}%
\right) .
\end{equation*}%
Since $u=\left[ 1,0,1,1\right] $ is the left eigenvector of $B$
corresponding to $\rho _{B}=2$, we have
\begin{eqnarray*}
uB_{1} &=&2\left[ 1,0,0,0\right] :=2u_{1}\mbox{ and }uB_{2}=2\left[ 0,0,1,1%
\right] :=2u_{2}, \\
u_{1}B_{2} &=&u_{2},\mbox{ }u_{2}B_{2}=u_{2},\mbox{
}u_{2}B_{1}=2u_{1}.
\end{eqnarray*}%
Therefore, $\mathcal{U}=\left\{ u_{1},u_{2}\right\} $ is the vector
space generated by vectors of the form $\left\{ JB_{J}:J\in
Y_{n}\mbox{ for }n\in \mathbb{N}\right\} $. Let $k=\dim
\mathcal{U}=2<4,$ and set
\begin{equation*}
L=\left(
\begin{array}{cccc}
1 & 0 & 0 & 0 \\
0 & 0 & 1 & 1%
\end{array}%
\right)
\end{equation*}%
be the $2\times 4$ matrix with rows forming a basis of $\mathcal{U}.$
Then
\begin{equation*}
LB_{1}=\left(
\begin{array}{cc}
0 & 0 \\
2 & 0%
\end{array}%
\right) \left(
\begin{array}{cccc}
1 & 0 & 0 & 0 \\
0 & 0 & 1 & 1%
\end{array}%
\right):= \widetilde{B}_1L \mbox{,}
\end{equation*}
\begin{equation*}
LB_{2}=\left(
\begin{array}{cc}
0 & 1 \\
0 & 1%
\end{array}%
\right) \left(
\begin{array}{cccc}
1 & 0 & 0 & 0 \\
0 & 0 & 1 & 1%
\end{array}%
\right) := \widetilde{B}_2L .
\end{equation*}%
Since $\dim \mathcal{U} = 2 = \#\mathcal{A}(Y)$, this implies $\left( l,\left\{ \widetilde{B}%
_{i}\right\} ,r\right) $ is indeed a presentation of the reduced
 module of the sofic measure $\nu $ on $Y$, where $l, r$ are
the left and right eigenvectors of $B$ corresponding to the maximal
eigenvalue $\rho_{\widetilde{B}} = \rho_B.$ Moreover, rank$\left( \widetilde{B}_{1}\right) =1$ and rank$\left( \widetilde{%
B}_{2}\right) =1$.
\end{example}

\begin{example}[\protect\cite{McMullen1984}]
\label{Ex: 1}Consider $\left( \mathbf{m},\mathbf{n}\right) =\left(
3,2\right) $. Let the adjacent matrix $A\in \mathbb{R}^{6\times 6}$
defining $Z$ as
\begin{equation*}
A=\left[
\begin{array}{cccccc}
1 & 0 & 0 & 1 & 1 & 0 \\
0 & 0 & 0 & 0 & 0 & 0 \\
0 & 0 & 0 & 0 & 0 & 0 \\
1 & 0 & 0 & 1 & 1 & 0 \\
1 & 0 & 0 & 1 & 1 & 0 \\
0 & 0 & 0 & 0 & 0 & 0%
\end{array}%
\right] .
\end{equation*}%
The permutation $P_{\left( 3,2\right) }$ induced from $\tau _{\left(
3,2\right) }$ (see Section 1) is
\begin{equation*}
P_{\left( 3,2\right) }\mathbb{=}\left[
\begin{array}{cccccc}
1 & 0 & 0 & 0 & 0 & 0 \\
0 & 0 & 1 & 0 & 0 & 0 \\
0 & 0 & 0 & 0 & 1 & 0 \\
0 & 1 & 0 & 0 & 0 & 0 \\
0 & 0 & 0 & 1 & 0 & 0 \\
0 & 0 & 0 & 0 & 0 & 1%
\end{array}%
\right] .
\end{equation*}%
Denote $P=P_{\left( 3,2\right) }$. Then the potential function $N$
extracted from $B=\left( B^{(i,j)}\right) _{i,j=1}^{2}$ by
(\ref{27}) is as follows.
\begin{equation*}
\left[
\begin{array}{cc}
\left( PAP^{-1}\right) ^{\left( 1,1\right) } & \left(
PAP^{-1}\right)
^{\left( 1,2\right) } \\
\left( PAP^{-1}\right) ^{\left( 2,1\right) } & \left(
PAP^{-1}\right)
^{\left( 2,2\right) }%
\end{array}%
\right] =\left[
\begin{array}{cccccc}
1 & 0 & 1 & 0 & 1 & 0 \\
0 & 0 & 0 & 0 & 0 & 0 \\
1 & 0 & 1 & 0 & 1 & 0 \\
0 & 0 & 0 & 0 & 0 & 0 \\
1 & 0 & 1 & 0 & 1 & 0 \\
0 & 0 & 0 & 0 & 0 & 0%
\end{array}%
\right] \allowbreak \mbox{.}
\end{equation*}%
One can easily check that $N=\left( N_{ij}\right)
_{i,j=1}^{2}=\left( B^{(i,j)}\right) _{i,j=1}^{2}$
\begin{eqnarray*}
N_{11} &=&\left(
\begin{array}{ccc}
1 & 0 & 1 \\
0 & 0 & 0 \\
1 & 0 & 1%
\end{array}%
\right) \mbox{, }N_{12}=\left(
\begin{array}{ccc}
0 & 1 & 0 \\
0 & 0 & 0 \\
0 & 1 & 0%
\end{array}%
\right) \mbox{, } \\
N_{21} &=&\left(
\begin{array}{ccc}
0 & 0 & 0 \\
1 & 0 & 1 \\
0 & 0 & 0%
\end{array}%
\right) \mbox{, }N_{22}=\left(
\begin{array}{ccc}
0 & 0 & 0 \\
0 & 1 & 0 \\
0 & 0 & 0%
\end{array}%
\right) \mbox{.}
\end{eqnarray*}%
is irreducible with
\begin{equation*}
V_{1}=\left(
\begin{array}{ccc}
1 & 0 & 1%
\end{array}%
\right) \mbox{ and }V_{2}=\left(
\begin{array}{ccc}
0 & 1 & 0%
\end{array}%
\right) \mbox{.}
\end{equation*}%
Since
\begin{equation*}
M= \left(
\begin{array}{cc}
2 & 2 \\
1 & 1%
\end{array}%
\right)
\end{equation*}%
Theorem \ref{Thm: 8} and Theorem \ref{Thm: 3} are thus applied to
show that the projection measure $\nu $ is a Markov measure which is
also the maximal measure of a full 2-shift. Finally, Theorem \ref{Thm: 6} is
also applied to show that
\begin{equation*}
\dim _{H}Z=\log_2 \rho _{M^{\alpha }}=\log_2 \left( 1+2^{\alpha }\right)
\mbox{, where }M^{\alpha }=\left(
\begin{array}{cc}
2^{\alpha } & 2^{\alpha } \\
1 & 1%
\end{array}%
\right) \mbox{.}
\end{equation*}
\end{example}

\section*{Acknowledgements}
The authors would like to thank D.-J.~Feng and L.-M.~Liao for their valuable suggestion during the preparation of this manuscript. Ban is supported in part by the National Science Council, ROC (Contract No NSC 99-2115-M-259-003), the National Center for Theoretical Sciences and the Center for Mathematics and Theoretical Physics at National Central University. Chang is grateful for the partial support of the National Science Council, ROC (Contract No NSC 100-2115-M-035-003).

\bibliographystyle{plainnat}

\bibliography{../../grece}

\end{document}